\documentclass[]{interact}
\usepackage[numbers,sort&compress]{natbib}
\bibpunct[, ]{[}{]}{,}{n}{,}{,}

\makeatletter
\def\NAT@def@citea{\def\@citea{\NAT@separator}}
\makeatother

\theoremstyle{plain}% Theorem-like structures
\newtheorem{theorem}{Theorem}[section]

\theoremstyle{definition}
\newtheorem{definition}[theorem]{Definition}

\theoremstyle{remark}

\usepackage{lineno,hyperref}
\usepackage{amsmath}
\usepackage{amssymb}
\usepackage{bm}
\usepackage{epsfig}
\usepackage{graphicx}
\usepackage[inline]{asymptote}
\usepackage{algorithm}
\usepackage{algpseudocode}
\usepackage{subfig}
\usepackage{paralist}
\usepackage{wasysym}
\usepackage{multirow}
\usepackage{float}

\makeatletter
\newsavebox{\sfe@box}
{\color@endgroup\egroup\subfloat[\sfe@caption]%
{\usebox{\sfe@box}}}
\makeatother

\modulolinenumbers[5]

\usepackage{color}
\usepackage{ulem}

\newcommand{\numcrit}{m}
\newcommand{\repset}{X_R}
\newcommand{\imrepset}{Z_R}
\newcommand{\R}{\mathbb{R}}
\newcommand{\N}{\mathbb{N}}
\newcommand{\lb}{l}
\newcommand{\ub}{u}
\newcommand{\lbs}{L}
\newcommand{\ubs}{U}
\newcommand{\uop}{\ubs^{\lb}}
\newcommand{\lop}{\lbs^{\ub}}
\newcommand{\ubsop}{\mathcal{\ubs}}
\newcommand{\lbsop}{\mathcal{\lbs}}
\newcommand{\uopkid}{\ubs^{\lb'}}

\newcommand{\commentout}[1]{}

\begin{document}

\articletype{ARTICLE TEMPLATE}

\title{An improved hyperboxing algorithm for calculating a Pareto front representation}
%\date{December 15, 2019}

\author{
\name{Kerstin D{\"a}chert\textsuperscript{a} 
and 
Katrin Teichert\textsuperscript{a}\thanks{CONTACT Katrin Teichert. Email: katrin.teichert@itwm.fraunhofer.de}}
\affil{\textsuperscript{a}Fraunhofer ITWM, Fraunhofer-Platz 1, Kaiserslautern, Germany}
}

\maketitle

\begin{abstract}
When solving optimization problems with multiple objective functions we are often faced with the situation that one or several objective functions are non-convex or that we can not easily show the convexity of all functions involved. In this case a general algorithm for computing a representation of the nondominated set is required. A suitable approach consists in a so-called \textit{hyperboxing algorithm} that is characterized by 
splitting the objective space into axis-parallel hyperrectangles. Thereby, only the property of nondominance is exploited for reducing the so-called search region. 
In the literature such an algorithm has already 
shown to provide a very good coverage of the Pareto front relative to the number of representation points calculated. However, the computational cost for the algorithm was prohibitive for problems with more than five objectives. 
In this paper, we present algorithmic advances that improve the performance of the algorithm and make it applicable to problems with up to nine objectives. 
We illustrate the performance gain and the quality of the representation for a set of test problems. We also apply the improved algorithm to a real world problem in the field of radiotherapy planning.
\end{abstract}


\section{Introduction}
Many real-world optimization problems contain more than one objective function. Hence, there is not one optimal objective function value but a set of nondominated points which build the so-called Pareto front. 
A decision maker usually wants to get an overview over this front to have an idea about the reachable alternatives to the given multiple criteria optimization (MCO) problem. Since the Pareto front can not be given in closed form, in general, a decision maker is typically satisfied with a finite representation, i.e.\ a finite number of nondominated points, or an approximation of the nondominated set. 
It is also possible to generate an inner and outer approximation and to sandwich the Pareto front in between. 
If the MCO problem is convex, weighted-sum scalarizations can be used to generate such an inner and outer approximation efficiently, see, e.g., \cite{BokFor12}. 
If the given problem contains non-convex functions or if the convexity of all functions can not be proven easily, a sandwich algorithm based on weighted-sum scalarizations is too restrictive in the sense that it cuts away feasible solutions on the one hand and is not able to reach every point on the other hand. 
Nevertheless, it is possible to generate an inner and an outer approximation of the nondominated set in the non-convex case, too. Firstly, the weighted-sum scalarization should be replaced by a scalarization that can, theoretically, generate every nondominated point. Secondly, we can only exploit the properties of nondominance when refining the inner and the outer approximation. This leads to the concept of \textit{hyperboxing algorithms}.

\begin{figure}[ht]
\centering
%\begin{subfloatenv}{ }
\subfloat[]
{{\includegraphics[width=.3\textwidth]{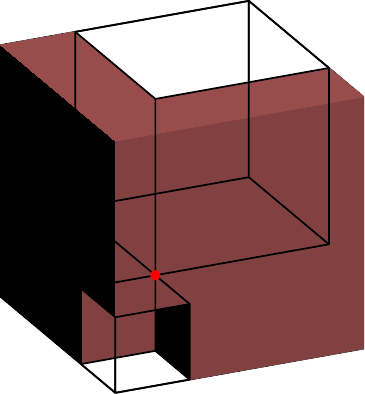} }}
%\end{subfloatenv}
\hspace{2cm}
%\begin{subfloatenv}{ }
\subfloat[]
{{\includegraphics[width=.3\textwidth]{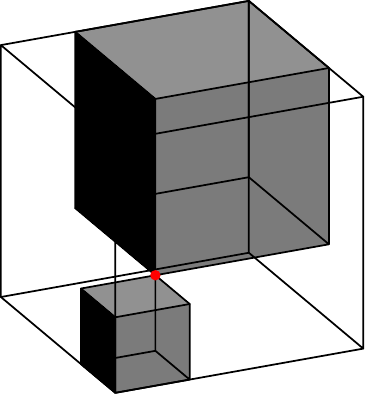} }}
%\end{subfloatenv}
\caption{Division of the objective space for a hyperboxing algorithm with respect to a single representation point $z$.
On the left, the remaining search space is depicted. On the right, the two boxes that are discarded are shown. 
The smaller grey box is known to be unattainable because $z$ is nondominated. The larger grey box is dominated by $z$ and thus cannot contain elements of the Pareto front either. 
}
\label{fig_hb_principle}
\end{figure}

Hyperboxing algorithms are a class of algorithms which construct a representation of the Pareto front, where exact Pareto solutions are computed in a deterministic, i.e.\ non-heuristic way. 
They are characterized in how - based on the representation points already found - they divide the objective space into parts that are kept and parts that can be discarded for the placement of additional representation points. These parts are the union of hyperboxes, hence the name. Figure \ref{fig_hb_principle} illustrates the division of the objective space for a hyperboxing algorithm under the assumption that the representation point $z$ is Pareto-optimal.

By cutting away unattainable regions from the objective space with each new representation point, an \textit{inner} and an \textit{outer approximation} of the Pareto front are defined. In case of minimization the lower border of the remaining search space represents an \textit{outer approximation}. 
Analogously, an \textit{inner approximation} for the Pareto front is defined by the upper border of the remaining search space. 
Hyperboxing algorithms may calculate and utilize either the inner approximation or the outer approximation or both for determining where to calculate the next representation point. The algorithm discussed in this paper uses both. 

Hyperboxing algorithms have been investigated by several authors. 
\cite{Sol91} describes the \textit{rectangle representation} of a bi-criteria nondominated set, which is equivalent to our hyperboxing approach in the 2D case. 
\cite{KlaTinWie02} present both a hyperboxing algorithm using the inner approximation and an algorithm using the outer approximation. 
\cite{Ser12} describes a hyperboxing algorithm that uses recursive splitting of the search boxes into disjoint parts. 
\cite{Tei14} 
uses a description of the search region by overlapping boxes each of which is defined by an inner and an outer knee point. When updating the set of inner and outer knee points, an algorithm similar to the 'redundancy elimination' strategy of \cite{KlaLacVan15} is used, i.e. the knee points are filtered at the end of each iteration to remove redundant ones. 
One improvement in this paper is to enhance the algorithm of \cite{Tei14}
with the 'redundancy avoidance' strategy of \cite{KlaLacVan15}.

The remainder of the paper is structured as follows. Section~\ref{sec:prelim} contains the notation and basic concepts. Section~\ref{sec:hypalg} explains the hyperboxing algorithm, followed by our algorithmic improvements in Section~\ref{sec:improve}. The latter are evaluated in Section~\ref{sec:comp} for academic test cases as well as a real-world problem from radiotherapy. In Section~\ref{sec:concl} we conclude and indicate directions for further research.  

\section{Preliminaries} \label{sec:prelim}

\noindent
\textit{Multi-criteria optimization} (MCO) seeks to solve a problem of the form 
\begin{equation} 
\begin{aligned}
\min_{x \in X} F(x) %\\ 
%x \in X
\end{aligned}
\end{equation}
where $x$ is the decision variable, $X$ is the feasible set, and
\begin{equation}
F=(f_1,..,f_{\numcrit}): X\rightarrow \mathbb{R}^{\numcrit}
\end{equation}
is a vector-valued function composed of real-valued objectives $f_1,..,f_{\numcrit}$. 
As the objectives are typically conflicting, there is no single solution to the MCO problem. Rather, we are interested in the set $X^{*}$ of best compromises whose elements are called \textit{Pareto efficient}: 

\begin{definition}
\label{def_Pareto_eff}
A feasible solution $x \in X$ is called \textit{Pareto efficient}, if there is no $\hat{x} \in X$ with $f_{i}(\hat{x}) \leq f_{i}(x)$ for all  $i \in \left\{1,..,{\numcrit}\right\}$ and $f_{j}(\hat{x}) < f_{j}(x)$ for some $j \in \left\{1,..,{\numcrit}\right\}$. 
We denote the set of efficient solutions by $X^{*}$. 
The set $Y^{*}=F(X^{*})$ is called the \textit{Pareto front} or the \textit{nondominated set}, its elements are called \textit{nondominated}.
\end{definition}  

For a continuous MCO problem, the set $X^{*}$ is generally infinite and cannot be calculated entirely. Therefore, MCO solution algorithms typically try to find a set of representative solutions $\repset \subseteq X$. According to \cite{Say00} the representativeness should be measured in the objective space: we want the set $F(\repset)$ to capture the Pareto front $Y^{*}$ as good as possible. We call $\imrepset:=F(\repset)$ a \textit{representation}. 

\begin{figure}[ht]
\centering
% \begin{asy}
%  include graph;
%  size(3cm);
%  pair z0 = (0.25, 0.8);
%  pair z1 = (0.4, 0.4);
%  pair z2 = (0.7, 0.23);
%  real alpha = 0.2;
%  pair s0 = (xpart(z0)-alpha, ypart(z0)-alpha);
%  pair s1 = (xpart(z1)-alpha, ypart(z1)-alpha);
%  pair s2 = (xpart(z2)-alpha, ypart(z2)-alpha);
%  path a0 = s0--(1., ypart(s0))--(1., 1.)--(xpart(s0), 1.)--cycle;
%  path a1 = s1--(1., ypart(s1))--(1., 1.)--(xpart(s1), 1.)--cycle;
%  path a2 = s2--(1., ypart(s2))--(1., 1.)--(xpart(s2), 1.)--cycle;
%  pen mygrey = rgb(0.9,0.9,0.9);
%  fill(a0, mygrey);
%  fill(a1, mygrey);
%  fill(a2, mygrey);
%  path f = (0.1, 0.9){right}..(0.25, 0.8)..(0.4,0.4)..(0.75, 0.2)..(0.9, 0.1);
%  draw(f);
%  dot(z0, black);
%  label("$z_{0}$",z0, NE);
%  dot(z1, black);
%  label("$z_{1}$",z1, NE);
%  dot(z2, black);
%  label("$z_{2}$",z2, NE);
%  label("$-\alpha e$", (0.,0.3), NE);
%  draw(z0--(xpart(z0)-alpha, ypart(z0)-alpha), EndArrow);
%  draw(z1--(xpart(z1)-alpha, ypart(z1)-alpha), EndArrow);
%  draw(z2--(xpart(z2)-alpha, ypart(z2)-alpha), EndArrow);
%  \end{asy}
\includegraphics[width=.3\textwidth]{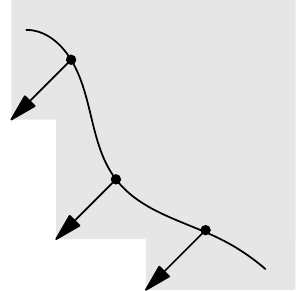} 
\caption{Example of an additive $\alpha$-approximation. The value of $\alpha$ determines the length of the arrows. }
\label{fig_alpha}
\end{figure}

\cite{Say00} discusses the aspects according to which the quality of a representation can be measured. \textit{Coverage} measures how far an element of $Y^{*}$ is maximally apart from the nearest representation point. \textit{Uniformity} measures how evenly spaced the representation points are. 
\textit{Cardinality} $\left|\imrepset \right|$ indicates the number of generated points, so implicitly it measures the computational effort for calculating $\imrepset$.
It is well known that attaining these three measures is a multi-objective problem with conflicting objectives 
itself, as, e.g., a good coverage typically implies a high cardinality. Hence, we have to decide how to balance these goals. 
The algorithm presented in this paper aims to achieve a maximally good coverage for a given cardinality (or, equivalently, tries to achieve a certain guaranteed coverage with  cardinality as low as possible), with coverage being measured based on the following definition (\cite*{BriFri10_1, PapYan00, Reu90}):

\begin{definition}
\label{def_additive_approx}
A representation set $\imrepset \subseteq F(X)$ is an \textit{additive $\alpha$-approximation} of $Y^{*}$ if for each $y \in Y^{*}$ there is $z \in \imrepset$ such that $z \leq y + \alpha e$, or, equivalently, if for each $y \in Y^{*}$ there is $z \in \imrepset$ such that
\begin{equation}
\label{eq_def_add_approx}
\max_{i = 1,..,\numcrit} (z_{i} - y_{i},0) =: d(y,z) \leq \alpha.
\end{equation}
We call the minimal $\alpha$ such that $\imrepset$ is an $\alpha$-approximation the \textit{approximation quality}.
\end{definition}
An illustration of the additive $\alpha$-approximation is given in Figure~\ref{fig_alpha}. The approximation quality as given in Definition~\ref{def_additive_approx} is a specific way of measuring the coverage in the sense of Say{\i}n \cite{Say00}, with equation (\ref{eq_def_add_approx}) defining the distance measure.

In order to find solutions to the MCO problem with conventional optimization algorithms, one has to aggregate the multiple objective functions into a single real-valued function. Such an aggregation function
\begin{equation}
\label{eq_scalarization}
\begin{aligned}
\sigma: Y\times \Lambda & \rightarrow \mathbb{R} \\
(F(x),\lambda) & \mapsto \sigma(F(x),\lambda)
\end{aligned}
\end{equation}
is called a \textit{scalarization}, with $\lambda \in \Lambda$ being the \textit{scalarization parameters}. 
The best known scalarization is probably the \textit{weighted sum scalarization} 
\begin{equation}
(F(x), w) \mapsto \sum^{\numcrit}_{i=1}w_{i}f_{i}(x).    
\end{equation}

Given a scalarization $\sigma$ an efficient (or at least weakly efficient) solution to the MCO can be found by minimizing $\sigma$ over the feasible set $X$ for certain fixed scalarization parameters $\lambda$. 
The weighted sum scalarization is easy to build and use but has the drawback that nondominated points not lying on the convex hull of the Pareto front can not be generated by it. Since we deal with arbitrary non-convex objective functions in this paper, we will use the \textit{Pascoletti-Serafini scalarization} introduced in \cite{PasSer84}. 
It is given by the mapping 
\begin{equation}
\begin{aligned}
(F(x),p,q) & \mapsto \alpha^{*} \\ 
\alpha^{*} & := \min \left\{\alpha \; | \ \ \exists \, \lambda \in \mathbb{R}^{\numcrit}_{+}: \ \ p + \alpha q = F(x) + \lambda \right\}. 
\end{aligned}
\end{equation}
Minimizing the Pascoletti-Serafini scalarization is equivalent to solving the following minimization problem, which we call the \textit{Pascoletti-Serafini problem} and denote by \textbf{PS}($p,q$):
\begin{equation}
\label{eq_ps_problem}
\begin{aligned}
\min_{\alpha, x, \lambda} \alpha \\
s.t. \ \ p + \alpha q & = F(x) + \lambda \\
(\alpha, x, \lambda) & \in \mathbb{R} \times X \times \mathbb{R}^{\numcrit}_{+}.
\end{aligned}
\end{equation}

\begin{figure}[ht]
\label{fig_pascoletti}
\centering
\includegraphics[width=.3\textwidth]{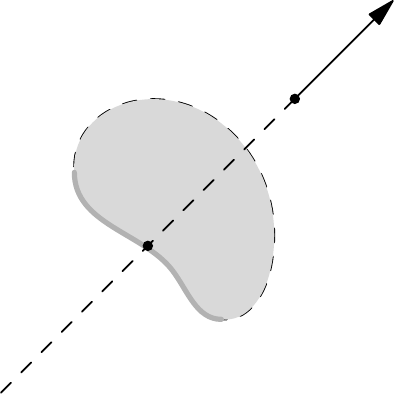} 
% \begin{asy}
%  include graph;
%  size(4cm);  
%  path y = (0.7,0.4)..(0.6,0.7)..(0.3,0.6)..(0.5,0.4)..(0.6,0.3)..cycle;
%  draw(y, dashed);
%  fill(y, rgb(0.85,0.85,0.85));
%  path ystar = (0.3,0.6){down}..(0.5,0.4)..{right}(0.6,0.3);
%  draw (ystar, linewidth(1.5)+rgb(0.7,0.7,0.7));
%  
%  pair p = (0.75, 0.75);
%  
%  dot(p, black);
%  label("$p$", p,SE);  
%  
%  real l = 0.2;
%  
%  draw(p--(xpart(p)+l, ypart(p)+l), EndArrow);
%  label("$q$", (xpart(p)+0.5*l, ypart(p)+0.5*l),SE);
%  
%  draw((xpart(p)-3*l, ypart(p)-3*l)--p, dashed);
%  dot((xpart(p)-0.3, ypart(p)-0.3), black);
%  label("$p+\alpha^{*}q$", (xpart(p)-0.3, ypart(p)-0.3),SE);
%  
%  label("$Y$", (0.45, 0.6),NE); 
%  \end{asy}
\caption{Illustration of the Pascoletti-Serafini problem $\textbf{PS}(p,q)$ where $\alpha^{*}$ is the optimal solution to \eqref{eq_ps_problem}.}  
\end{figure}

%\subsection{Concept of the algorithm}
%%%%%%%%%%%%%%%%%%%%%%%%%%%%%%%%%%%%%%%%%%%%%%%%%%%%%%%%%%%%%%%%%%%%%%%%%%%%%%%%%%%%%%%%%
%
% The hyperboxing algorithm
%
%%%%%%%%%%%%%%%%%%%%%%%%%%%%%%%%%%%%%%%%%%%%%%%%%%%%%%%%%%%%%%%%%%%%%%%%%%%%%%%%%%%%%%%%%

\section{The hyperboxing algorithm} \label{sec:hypalg}
\subsection{General concept}

\begin{figure}[ht]
\centering
\subfloat[]
{{\includegraphics[width=.27\textwidth]{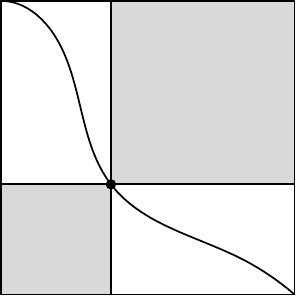} }}
\hspace{.5cm}
\subfloat[]
{{\includegraphics[width=.27\textwidth]{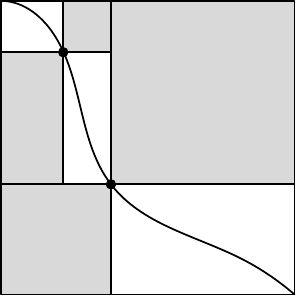} }}
\hspace{.5cm}
\subfloat[]
{{\includegraphics[width=.27\textwidth]{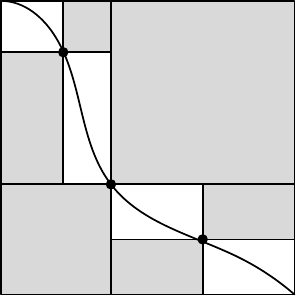} }}
\caption{Workflow of the hyperboxing algorithm for a bi-objective problem. The boxes enclosing the Pareto front are depicted in white, the excluded parts of the outcome space are depicted in grey. In (b), the box at the bottom right has clearly the largest minimal edge length and is therefore selected for further refinement in (c). }
\label{fig_hb_2d}
\end{figure}

In principle, the hyperboxing algorithm discussed in this paper works as follows. At the beginning, a  \textit{start box} $B_{0}\subseteq \mathbb{R}^{\numcrit}$ is defined that contains the part of the Pareto front one wants to approximate. 
In the first iteration, the Pascoletti-Serafini problem is solved where  parameter $p$ denotes the upper corner and $q$ the diagonal of $B_{0}$. 
Let $(\alpha_{0}, x_{0}, \lambda_{0})$ be the solution found, then 
$z_{0}=F(x_{0})$ is added to the representation set $\imrepset$. 
Now, using the property that $z_{0}$ is nondominated, 
$B_{0}$ can be decomposed into new smaller boxes
(how this is done will be detailed below). 
From those, the largest one $B_{1}$ is chosen, where the size of a box is 
defined as the smallest edge of the box. Details on why we choose this measure for the box size will also be given below.
Again, the upper corner of $B_{1}$ and its diagonal define a Pascoletti-Serafini problem, which is solved to obtain the new representative solution $x_{1}$ and representation point $z_{1}$, and so on. For a two dimensional Pareto front, the working principle of the algorithm is illustrated in Figure \ref{fig_hb_2d}.
A basic implementation is given in Algorithm~\ref{alg_ahb_general}.
The computationally expensive tasks are the update procedures \textit{newLowerBounds} and \textit{newUpperBounds} as well as finding the largest box from the potentially very large set $\mathcal{B}$. They are explained in more detail in Sections~\ref{subsec:updatebounds} and \ref{subsec:refine} and give rise to our algorithmical improvements which are presented in Section~\ref{sec:improve}. 

\begin{algorithm}[htbp]           
\caption{Basic hyperboxing algorithm}          
\label{alg_ahb_general} 
\begin{algorithmic}%[1]
\algblockdefx[Input]{Input}{EndInput}
    {\textbf{input:}}
    {}
\algblockdefx[Output]{Output}{EndOutput}
    {\textbf{output:}}
    {}
\Input
\begin{compactitem}
	\item a multi-objective problem with $\numcrit$ objectives
	\item an initial box $B_{0}=\left[\lb_0,\ub_0 \right]\subseteq \mathbb{R}^{\numcrit}$
	\item a target approximation quality $\epsilon > 0$
\end{compactitem}
\EndInput
\Output
\begin{compactitem}
	\item a set of (weakly) nondominated representation points $Z$
\end{compactitem}
\EndOutput

\algblock[Name]{Start}{End}
\Start
\State $Z \gets \emptyset$
\State $\mathcal{B} \gets B_{0}$
%\State $\lbs \gets \left\{b^{l}\right\}$, $\ubs \gets \left\{b^{u}\right\}$
\State $\lbs \gets \left\{\lb_0 \right\}$, $\ubs \gets \left\{\ub_0 \right\}$
\While {$\max_{B \in \mathcal{B}} \text{size}(B) > \epsilon$}
		\State Pick $B = [\lb,\ub]$ from $\mathcal{B}$ such that $\text{size}(B)$ is maximal
		\State solve $PS(\ub, \ub-\lb)$ and obtain solution $(\alpha, z, \lambda)$
		\State $s := z +\lambda$
		\State $Z \gets Z\cup \left\{z\right\} $
		\State $\lbs \gets newLowerBounds(\lbs, s)$ 
		\State $\ubs \gets newUpperBounds(\ubs, z)$
		\State $\mathcal{B} := \left\{B = \left[\lb,\ub\right] | \ \ \lb \in \lbs, \ub \in \ubs, \lb < \ub \right\}$
\EndWhile
\State \Return $Z$
\End
\end{algorithmic}
\end{algorithm}

%%%%%%%%%%%%%%%%%%%%%%%%%%%%%%%%%%%%%%%%%%%%%%%%%%%%%%%%%%%%%%%%%%%%%%%%%%%%%%%%%%%%%%%%%
\subsection{Relationship between boxes, inner and outer approximations and (local) lower and upper bounds}

\begin{figure}[htbp] 
  \centering
\subfloat[]
{{\includegraphics[width=.3\textwidth]{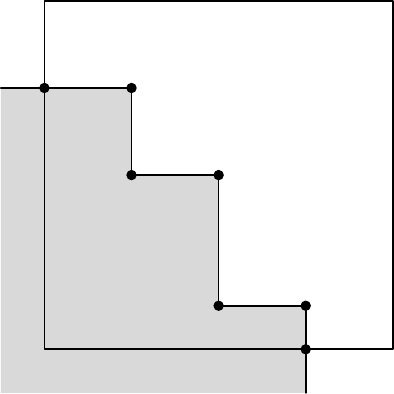} }}
\hspace{1cm}
\subfloat[]
{{\includegraphics[width=.3\textwidth]{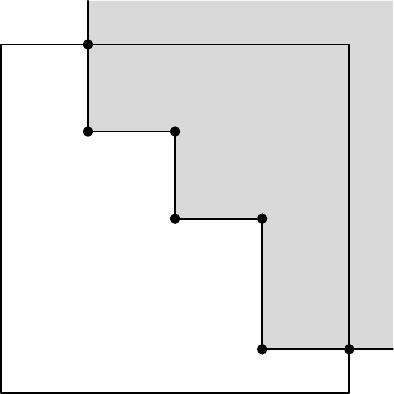} }}
	\caption{(a) The lower bounds $\lbs = \left\{\lb_{1},..,\lb_{4}\right\}$ and (b) the upper bounds $\ubs = \left\{\ub_{1},..,\ub_{4}\right\}$ for a set 
	$\imrepset = \left\{z_{1},z_{2},z_{3}\right\}$ with respect to some start box $B_{0}$. 
	The connecting lines in (a) build an inner approximation, the lines in (b) an outer approximation. 
	Note that the update of the start box with respect to $\imrepset$ is not depicted. }
	\label{fig_lbub}
\end{figure}

In the convex case when a sandwich algorithm is used to enclose the Pareto front in between an inner and an outer approximation, these approximations have to be constructed explicitly from the previous inner and outer approximation, the current nondominated point and the scalarization parameters. 
In the non-convex case the situation becomes easier in the sense that less information has to be stored. In fact, besides the nondominated points, it is enough to save and update the lower and upper corner points of the boxes which we refer to as (local) lower and upper bounds \cite{KlaLacVan15} in the following. 
\begin{definition}[Lower and upper bounds]%\cite{LegLeGCot10}
\label{def:bounds}
Let $B_{0} \subset \mathbb{R}^{\numcrit}$ be a box and $Z\subseteq B$ a finite set. 
Then the minimal elements of $B_{0}\setminus (Z - \mathbb{R}^{\numcrit}_{+})$ are called the \textit{lower bounds} (with respect to $B_{0}$) and denoted by $\lbs$. 
The maximal elements of $B_{0} \setminus (Z + \mathbb{R}^{\numcrit}_{+})$ are called the \textit{upper bounds} and denoted by $\ubs$. 
\end{definition}
In the literature there is also the alternative notion of \textit{knee points}
\cite{LegLeGCot10}. 
Figure~\ref{fig_lbub} shows the set of local lower and upper bounds with respect to a given set of nondominated points in the bicriteria case. 
The inner and outer approximation can be easily reconstructed from these points, however, this is not required for the course of the algorithm.
For the latter, it is enough to know the current set of local lower and upper bounds as well as the information which pairs of lower and upper bounds build a box. 

%%%%%%%%%%%%%%%%%%%%%%%%%%%%%%%%%%%%%%%%%%%%%%%%%%%%%%%%%%%%%%%%%%%%%%%%%%%%%%%%%%%%%%%%%
\subsection{Updating the lower and upper bounds} \label{subsec:updatebounds}

In the following, we discuss the update of the upper bounds. The update of the lower bounds works in a similar fashion, hence we only highlight whenever differences occur. 

Figure~\ref{fig_hb_2d} illustrates the update in the bi-objective case. Each time we find a new representation point $z$, we split the box containing this point into two new boxes. This means that the upper bound $\ub$, which is strictly greater than $z$, 
is replaced by new upper bounds $\ub_1$ and $\ub_2$.
In this process, $\ub_1$ inherits the values of the former upper bound $\ub$ in the second component, whereas $\ub_2$ inherits the values of $\ub$ in the first component.
The component that is not inherited from $\ub$ is filled with the value of the new point $z$ in the respective component, i.e.\ the first component of $\ub_1$ equals the first component of $z$ and the second component of $\ub_2$ equals the second component of $z$.

Analogously, if a lower bound $\lb$ is strictly smaller than $z$ in every component, we replace it by two new lower bounds $\lb_1$ and $\lb_2$ . 
Note that for calculating the components of the lower bounds we can use the point 
$s = z + \lambda$ coming from the Pascoletti-Serafini scalarization instead of $z$.
In case that the diagonal of the box does not intersect with the Pareto front, the parameter $\lambda$ is positive, i.e., using $s$ creates tighter bounds than using $z$ (otherwise $s=z$ holds.)  

In higher dimensions the same principle applies concerning the update of a bound, i.e.\ a new bound inherits all but one component from the bound it was created from, and the new point ($z$ or $s$, respectively) determines exactly one component of the new bound. 
However, updating the bounds becomes much more complicated,
since a new point can be situated in more than one box, and therefore require the replacement of multiple lower and upper bounds. 
 
If we decompose every bound which is 'affected' by a the new representation point into $\numcrit$ new bounds, we obtain redundant bounds \cite{DaeKla15, KlaLacVan15}. 
A viable approach to deal with these redundant bounds is to filter them out afterwards by pairwise comparisons \cite{Tei14}.
An improvement consists in avoiding the creation of redundant bounds beforehand. How to do so is described in Section~\ref{sec:improve} below.

\subsection{Refinement of the approximation} \label{subsec:refine}

Finally, in every iteration we want to improve the approximation as best as possible. Practically, we want to select that pair of a local lower and upper bound that builds the box with the largest size. In principle, various choices for measuring the box size are feasible, e.g. the volume.

\begin{figure}[]
\centering
%\begin{asy}
%  include graph;
%  size(4cm);
%  pair z0 = (0.5, 0.8);
%  pair z1 = (0.8, 0.4);
%  real l2 = 0.4; 
%  real l1 = 0.3;
%  pair e0 = (xpart(z0)-l1, ypart(z0)-l1);
%  pair e1 = (xpart(z1)-l1, ypart(z1)-l1);
%  path f = (0.4, 0.9){right}..z0{down}--(0.51, 0.8)..(0.52, 0.6)..(0.53, 0.5)..z1..(0.9, 0.3);
%  path b = z0--(xpart(z0), ypart(z1))--z1--(xpart(z1), ypart(z0))--cycle;
%  path gb = e1--(1., ypart(e1))--(1., 1.)--(xpart(e1), 1.)--cycle;
%  path gb2 = e0--(1., ypart(e0))--(1., 1.)--(xpart(e0), 1.)--cycle;
%  pen mygrey = rgb(0.9,0.9,0.9);
%  fill(gb, mygrey);
%  fill(gb2, mygrey);
%  draw(b);
%  draw(f);
%  dot(z0, black);
%  label("$z_{0}$",z0, NE);
%  dot(z1, black);
%  label("$z_{1}$",z1, NE);
%  real l2 = 0.4; 
%  real l1 = 0.3;
%  draw(z0--e0, EndArrow);
%  draw(z1--e1, EndArrow);
%  label("$l_1$", (0.53,0.28), NE);
%  label("$l_2$", (0.4,0.5), NE);
%  label("$-l_1 e$", (0.57,0.1), NE);
%\end{asy}
{{\includegraphics[width=.3\textwidth]{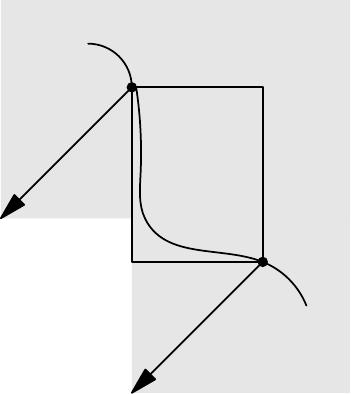} }}
\caption{Motivation for using the smallest edge length as measurement for the box size. If we only consider the part of the Pareto front within the box, then the smallest edge length $l_{1}$ of the box is an upper bound for the approximation quality according to Definition~\ref{def_additive_approx}. There is a representation point $z_1$ such that $z_1 \leq y + l_{1} e$ for all nondominated $y$ within the box no matter how the (unknown) Pareto front is shaped.}
\label{fig_min_side_length_alpha}
\end{figure}

A natural choice for measuring the size of a box is its smallest edge length. With this choice, the size of the largest box is an upper bound for the approximation quality according to Definition~\ref{def_additive_approx}. Figure~\ref{fig_min_side_length_alpha} illustrates this correspondence between box size and approximation quality for the bi-criteria case. As the smallest side length of each box provides an upper bound for the approximation quality for the part of the Pareto front inside the box, computing the maximum over the smallest edge lengths of all boxes yields an upper bound for the global approximation quality. For an exact formulation and proof of the correspondence between this choice of measuring the box size and the additive approximation quality, see \cite{Tei14} (Corollary 5.1.9).

Remind that, in order to define a box, the lower bound must be strictly smaller than the upper bound in every component. We call a bound satisfying this property \textit{opposing} in the following. 
Again, the bicriteria case is much simpler than the general multi-objective case. In the bicriteria case, for every bound there is a unique opposing bound. However, for three and more criteria, 
more than one upper bound can be strictly greater than a certain lower bound, and vice versa. 
As an example consider Figure~\ref{fig_hb_principle} with three lower and three upper bounds. For every upper bound there are two lower bounds that are strictly smaller in every component. 
Hence, for more than two criteria, for every  bound we obtain not only one but a list of opposing bounds. 

Naively, we can store the lower bounds $\lbs$ as a list and the upper bounds $\ubs$ as a second list. Then in each iteration we identify the largest box by iterating over the nested lists. This approach is used in \cite{Tei14}. How to improve it is discussed in Section~\ref{sec:improve} below. 

%%%%%%%%%%%%%%%%%%%%%%%%%%%%%%%%%%%%%%%%%%%%%%%%%%%%%%%%%%%%%%%%%%%%%%%%%%%%%%%%%%%%%%%%%
\section{Algorithmic improvements} \label{sec:improve}
%\subsubsection{Increasing algorithmic efficiency by storing opposing bounds}

%%%%%%%%%%%%%%%%%%%%%%%%%%%%%%%%%%%%%%%%%%%%%%%%%%%%%%%%%%%%%%%%%%%%%%%%%%%%%%%%%%%%%%%%%
\subsection{Efficient update of the bounds using defining points}
\label{sec_eff_update}

As explained above, a bound is defined in each component by the respective component of a point~$z$ (or $s$). 
(Only bounds at the boundary might be   
defined by the components of the start box.) 
\cite{KlaLacVan15} make use of this fact and store these so-called \textit{defining points}. 
Based on them they formulate a criterion with the help of which redundant bounds can be identified before their creation. More precisely, for a certain upper bound $\ub$ greater than the current representation point $z$ and for a certain component $k$, they consider the value of the $k$th component of all but one of the defining points, leaving out the $k$th defining point. Then, they compute the maximum of these $\numcrit-1$ values. The new bound $\ub_{k}$, which would have the $k$th component component set to $z_{k}$, is then created if and only if $z_{k}$ is strictly greater than this computed maximum. 

In \cite{KlaLacVan15} only the upper bounds and their update are required. 
However, their idea can be applied to the update of the lower bounds in a straightforward way. 
Instead of building the maximum, we compute the minimum of the defining points and then test whether $s_{k}$ is strictly smaller than this computed minimum.

%%%%%%%%%%%%%%%%%%%%%%%%%%%%%%%%%%%%%%%%%%%%%%%%%%%%%%%%%%%%%%%%%%%%%%%%%%%%%%%%%%%%%%%%%
\subsection{Storing opposing bounds}
\label{subsec_opposing_point_sets}

Another important computational improvement consists in 
storing \textit{opposing bounds} as follows: 
 For each lower bound $\lb \in \lbs$, denote
\begin{equation}
    \uop = \left\{\ub \in \ubs \, | \, \lb < \ub \right\}
\end{equation}
the set of opposing upper bounds; and for each upper bound $\ub \in \ubs$, denote
\begin{equation}
    \lop = \left\{\lb \in \lbs \, | \, \lb < \ub \right\}
\end{equation}
the set of opposing lower bounds. 
Using this approach, whenever we update either $\lbs$ or $\ubs$, we also have to update 
$\ubsop = \left\{\uop \, | \,  \lb \in \lbs \right\}$ 
and 
$\lbsop = \left\{\lop \, | \,  \ub \in \ubs \right\}$. 

There are different possibilities in which order the sets can be updated. We propose to first update $\lbs$ and then (temporarily) $\ubsop$ and $\lbsop$, followed by the update of $\ubs$ and then again $\lbsop$ and $\ubsop$, however now in reverse order. Of course, the two main steps are interchangeable, however, the intermediate update of $\ubsop$ and $\lbsop$ in the first step is required for not losing the correct relationships among the opposing bounds. 

In the following, we explain in more detail 
the first part of this update, i.e. the temporary update of $\ubsop$ and $\lbsop$. The second step after the update of $\ubs$ then works in a completely analogous way. 

Assume that $\lbs$ has been updated. Regarding the update of $\ubsop$ %(step 2) 
note that the components of the lower bounds are monotonically increasing over the course of the iterations. 
Namely, if a lower bound $\lb$ is strictly smaller than a newly found point $s$, in the update process it is replaced by a set of child lower bounds $\lb'_{1},...,\lb'_{n}, n \leq \numcrit,$ created from $\lb$ and $s$ as described in Section~\ref{subsec:updatebounds} and in more detail in \cite{KlaLacVan15}. 
Since $\lb' \geq \lb$ component-wise, we have $\uopkid \subseteq \uop$ for each child $\lb'$ of $\lb$. 
Therefore, we can calculate $\uopkid$ for a new node $\lb'$ by filtering over $\uop$, i.e. removing every $\ub$ that does not fulfill $\lb' < \ub$. 

Regarding the subsequent update of $\lbsop$, %(step 3), 
note that $\uop$ is exactly the set of upper bounds $\ub$ for which $\lb \in \lop$ holds. Also as $\lb \leq \lb'$ for all children, $\lb' \in \lop$ cannot hold if not $\lb \in \lop$ before. 
Therefore, we only have to update the sets $\lop$ for $\ub \in \uop$. 
To this end, for all $\ub \in \uop$ we remove $\lb$ from $\lop$ and insert all children $\lb'$ of $\lb$ with $\lb' < \ub$. 

In the second step we update $\ubs$ with respect to $z$ and then again $\lbsop$ and $\ubsop$.
After this update, we can easily determine the largest box. 
To this end, we first find for each $\lb \in \lbs$ the $\ub^{*} \in \uop$ such that the size of the corresponding box is maximized, and then we pick the largest box 
among all these boxes. 

This algorithmic approach can be improved even more by observing that we do not need to store all opposing bounds. 
For a lower bound $\lb$ we only have to store an opposing upper bound $\ub$ if the box size of $B = [\lb,\ub]$ is larger than the target approximation quality $\epsilon$. 
Therefore, when we filter for $\lb' < \ub$ in the creation of $\uopkid$ and $\lop$ as described above, %(steps 2 and 3), 
we can directly exclude boxes $B = [\lb',\ub]$ whose minimal edge length is smaller than $\epsilon$. 
The same can be done after the update of $\ubs$.
This reduces the sizes of the sets $\uop$ and $\lop$ considerably and becomes especially powerful in later iterations when nearing the approximation quality (and when speedup does matter most). 
%See the appendix for a description of the algorithm in pseudo code.

%%%%%%%%%%%%%%%%%%%%%%%%%%%%%%%%%%%%%%%%%%%%%%%%%%%%%%%%%%%%%%%%%%%%%%%%%%%%%%%%%%%%%%%%%
%
% Comp study
%
%%%%%%%%%%%%%%%%%%%%%%%%%%%%%%%%%%%%%%%%%%%%%%%%%%%%%%%%%%%%%%%%%%%%%%%%%%%%%%%%%%%%%%%%%

\section{Computational study} \label{sec:comp}

\subsection{Setup}

The hyperboxing algorithm is implemented in C++. The interface of this implementation allows to iteratively add points to the representation and get the resulting largest box, which then defines the parameters of the next Pascoletti-Serafini problem to solve. The Pascoletti-Serafini problems themselves are solved within Matlab, either analytically (for the 3D sphere and the hyperellipsoid example) or using the Knitro solver \cite{ByrNocWal06}. 
All plots are generated in Matlab. 

As initial lower bound we use the \textit{ideal point} which is a tight lower bound on the nondominated set. It can be computed rather easily by the (component-wise) determination of the  minimum values of the individual objective functions. As initial upper bound we use the \textit{nadir point} which is a tight upper bound on the nondominated set. 
Note that the latter is often not available beforehand since it requires knowledge of the nondominated set. Instead, any other meaningful initial upper bound can be selected. 

\subsection{Test Problems}

As test problems we consider five test cases for three objectives. 
Two of them are directly formulated for any number of objectives and serve as the benchmark problems for our numerical study with four to nine objectives. 

\subsubsection{(Hyper-)Sphere and (Hyper-)Ellipsoid}
\noindent 
The formulation of both test problems is  
\begin{equation}
\label{testproblem_hyper}
\begin{aligned}
\min &
\left( 
\begin{array}{c}
      x_{1} \\ 
      x_{2} \\
      \vdots \\
      x_{\numcrit}
\end{array}
\right)
\\
s.t \ \ \sum^{\numcrit}_{i = 1} \left(\frac{x_{i}}{a_{i}}\right)^{2} & \leq 1, \ \ \text{with}\ \ a_i > 0 \ \ \forall i. \\
\end{aligned}
\end{equation}
For the (hyper-)sphere problem all  parameters $a_1, a_2, \dots, a_{\numcrit}$ are set to $1$.
For the (hyper-)ellipsoid we set $a_1=\numcrit$ and $a_2, \dots, a_{\numcrit} = 1$.

For both problems the feasible set equals the feasible outcome set and is either a unit ball or an ellipsoid centered at the origin. 
The Pareto front is that part of the boundary that lies in the 'lower left' part, i.e.\ with all components smaller or equal to zero. 
The ideal point is, thus, 
$(-1,-1,\dots,-1)^{\top}$, the nadir point 
$(0,0,\dots,0)^{\top}$. 

Moreover, the Pareto front is $\mathbb{R}^{\numcrit}_+$-convex, and its Pascoletti-Serafini scalarization problems can be solved analytically, making this a test case that is simple to solve.
Namely, if $A$ denotes the diagonal matrix with diagonal entries $a_{i}$, 
the boundary of the feasible set in the objective space can be described by all $z \in \R^{\numcrit}$ for which $z^{\top} A z = 1$ holds.
On the other hand, any feasible solution of the Pascoletti-Serafini problem \textbf{PS}($p,q$) must satisfy the constraint $p + \alpha q = F(x) + \lambda$ from \eqref{eq_ps_problem}. Since \eqref{testproblem_hyper} is convex, $\lambda = 0$ and, thus, $z = F(x) = p + \alpha q$ holds. Hence, we can find the minimal $\alpha^{*}$ of the Pascoletti-Serafini problem \textbf{PS}($p,q$) by solving the quadratic equation
\begin{equation}
\begin{aligned}
0 & = (p + \alpha q)^{\top}A(p + \alpha q) - 1 \\
  & = (q^{\top}Aq)\alpha^{2} + (2p^{\top}Aq)\alpha + (p^{\top}Ap - 1).
\end{aligned}
\end{equation}

\subsubsection{Non-convex connected front}
\noindent Our third test case (\ref{testproblem_eichfelder}) has a connected but $\mathbb{R}^3$-non-convex Pareto front. 
It has been used in \cite{Eic06} and \cite{gourion10} and is a slight variation of a test problem from \cite{KimWec06}. 

\begin{equation}
\label{testproblem_eichfelder}
\begin{aligned}
\min 
\left( 
\begin{array}{c}
      -x_{1} \\ -x_{2} \\ -(x_{3})^{2} 
\end{array}
\right) &
\\
s.t. \ \ -\cos{(x_{1})} - \exp{(-x_{2})} + x_{3} & \leq 0 \\
0 & \leq x_{1} \leq \pi \\
0 & \leq x_{2} \\
1.2 & \leq x_{3} \\
\end{aligned}
\end{equation}
According to \cite{gourion10} its efficient set is given by 
\begin{align*}
\{ 
x \in \R^3 : \; &0 \leq x_{1} \leq \arccos(0.2), 
0 \leq x_{2} \leq -\ln(1.2-\cos(x_{1})),\\ 
&1.2 \leq x_{3} \leq \cos(x_{1})+\exp(-x_{2}) 
\}.
\end{align*}
The ideal point 
%$z^I$ %= -(\arccos(0.2),-\ln(0.2),4)^{\top} \approx 
is approximately 
$(-1.37,-1.61,-4)^{\top}$, the nadir point 
%z^N = 
$(0,0,-1.44)^{\top}$.

%\begin{figure}[H]
%\begin{center}
%\includegraphics[width=0.5\linewidth]{Eichfelder_paretofront.png}
%\end{center}
%\caption{Pareto front of test problem \eqref{testproblem_eichfelder}}
%\label{fig_eichfelder}
%\end{figure}

\subsubsection{Comet}
\noindent This test problem (\ref{testproblem_deb1}) stems from \cite{DebThiLauZit01}. The Pareto front resembles a comet. %, as illustrated in Figure~\ref{fig_deb1}.

\begin{equation}
\label{testproblem_deb1}
\begin{aligned}
\min &
\left( 
\begin{array}{c}
      (1 + x_{3})(x_{1}^{3}x_{2}^{2} - 10 x_{1} - 4 x_{2}) \\ 
      (1 + x_{3})(x_{1}^{3}x_{2}^{2} - 10 x_{1} + 4 x_{2}) \\ 
      3(1 + x_{3})x_{1}^{2} 
\end{array}
\right)
\\
s.t \ \ 1 & \leq x_{1} \leq 3.5 \\
-2 & \leq x_{2} \leq 2 \\
0 & \leq x_{3} \leq 1 \\
\end{aligned}
\end{equation}
In~\cite{DebThiLauZit01} it is stated that the set
$$
S := \{ x \in \R^3 : 
1 \leq {x}_{1} \leq 3.5, \; -2 \leq {x}_{2} {x}_{1}^3 \leq 2, \; {x}_{3}=0\} 
$$ 
represents the efficient set. This is, however, not correct. 
Consider, e.g., $\bar{x}:=(3.5,0,0)^{\top}$. 
Then $\bar{x} \in S$ and $
f(\bar{x})=(-35,-35,36.75)^{\top}$. 
Let $\tilde{x}:=(2,0,1)^{\top}$. 
This point is feasible for \eqref{testproblem_deb1} and $\tilde{x} \notin S$ holds. 
However, $f(\tilde{x})=(-40,-40,24)^{\top} \leqq f(\bar{x})$. Hence, $S$ does not represent the efficient set of \eqref{testproblem_deb1}. Without proof, but confirmed by the numerical results, we claim that the 
nondominated set of \eqref{testproblem_deb1} is given by $f(S')$ with 
\begin{align} \label{chap5:comet:xe}
S' := & \{ x \in \R^3 : 
1 \leq {x}_{1} \leq 3.5, \; -2 \leq {x}_{2} {x}_{1}^3 \leq 2, \; {x}_{3}=1
\} \\ \notag
& \cup
\{ x \in \R^3 : 
{x}_{1} =1, \; -2 \leq {x}_{2} \leq 2, \; 0 \leq {x}_{3} \leq 1
\}.
\end{align}
Note that the same shape of the nondominated set has been obtained in the numerical study of~\cite{Eic06}. 
The ideal point is 
%$z^I \approx 
approximately $(-70.19,-70.19,3)^{\top}$, 
%the vector of individual maxima is $z^M=(289,289,73.5)^{\top}$. 
%The local nadir point is given by $(2,2,73.5)^{\top}$,
the nadir point is 
%$z^N=
$(4,4,73.5)^{\top}$.

%\begin{figure}[H]
%\begin{center}
%\includegraphics[width=0.5\linewidth]{Deb1_paretoFront.png}
%\end{center}
%\caption{Pareto front of test problem \eqref{testproblem_deb1}}
%\label{fig_deb1}
%\end{figure}

\subsubsection{Patched Pareto front}
This test problem is a modification of problem DTLZ7 \cite{DebThiLauZit01}. 
The problem is designed for an arbitrary number of objectives and $n=m-1+k$ variables, where $k \in \N$ is a parameter. The authors suggest $k=20$.
The resulting problem formulation in the tricriteria case, which is also considered in \cite{Eic06}, is 
\begin{equation}\label{chap5:disconnect}
\begin{array}{ll} 
\min & \left(
\begin{array}{c}
x_{1}\\
x_{2}\\
g(x) \cdot \left( 
3- \displaystyle \sum_{i=1}^2 \left( \frac{x_{i}}{g(x)} \left( 1+\sin(3 \pi x_{i} \right) \right) 
\right)
\end{array}
\right)
\\
\mbox{s.t.} & \\[-0.3cm]
& x_{i} \in [0,1] \quad \forall \, i=1,\dots,2+k,\\
& x \in \R^{2+k}
\end{array}
\end{equation}
with 
$$
g(x)=2+\frac{9}{k} \displaystyle \sum_{i=3}^{2+k} x_{i}.
$$
For $k=20$, a nonlinear solver experiences problems when solving the scalarized problems. 
As our focus does not lie on improving existing nonlinear single objective solvers, we consider the modified problem
\begin{equation}%\label{chap5:disconnectsimpl}
\label{testproblem_deb2}
\begin{array}{ll} 
\min & \left(
\begin{array}{c}
x_{1}\\
x_{2}\\
6- \displaystyle \sum_{i=1}^2 \left( x_{i} \left( 1+\sin(3 \pi x_{i} \right) \right) 
\end{array}
\right)
\\
\mbox{s.t.} & \\[-0.3cm]
& x_{i} \in [0,1] \quad \forall \, i=1,2,\\
& x \in \R^{2}
\end{array}
\end{equation}
with only two variables. 
Note that \eqref{testproblem_deb2} is not a special case of~\eqref{chap5:disconnect}, as $g(x)$ is not defined for $k=0$. However, the nondominated set of~\eqref{testproblem_deb2} is the same as of~\eqref{chap5:disconnect} for arbitrary $k \in \N_{0}$, which can be easily seen. 
The third objective of~\eqref{chap5:disconnect} can be reformulated as
$$
f_{3}(x) = 3 g(x) - \displaystyle \sum_{i=1}^2 x_{i} \left( 1+\sin(3 \pi x_{i} \right) ).
$$
As $g(x)$ does not depend on $x_{1}$ and $x_{2}$, 
and, at the same time, the variables $x_{3}, \dots,x_{k+2} \in [0,1]$ only occur in $g(x)$, which only occurs in the third objective function, we can eliminate the variables 
$x_{3},\dots,x_{k+2}$ from \eqref{chap5:disconnect}. 
Indeed, for any nondominated point of \eqref{chap5:disconnect}, 
%, $f_{3}(x)$ becomes minimal for 
$x_{3}=\dots=x_{k+2}=0$ must hold.
%(see also \cite{deb01b}, \citet{deb05}, \cite{Eic06}).
This implies $g(x)=2$. 
As $f(x)=(x_{1},x_{2},f_{3}(x))^{\top}$ only depends on $x_{1}$ and $x_{2}$, the nondominated sets of~\eqref{chap5:disconnect} and 
\eqref{testproblem_deb2} are the same. 
It consists of four disconnected parts. 
The ideal point is
%$z^I=
$(0,0,2.61)^{\top}$,
the nadir point is 
%$z^N = 
$(0.86,0.86,6)^{\top}$.

\commentout{
\noindent This test problem (\ref{testproblem_deb2}) is a variation of a test case from \cite{DebThiLauZit01}. The Pareto front consists of four disconnected patches. 
%and is illustrated in Figure~\ref{fig_deb2}.

\begin{equation}
\label{testproblem_deb2}
\begin{aligned}
\min 
\left( 
\begin{array}{c}
      x_{1} \\ 
      x_{2} \\ 
      6 - \sum_{i=1}^{2} x_i (1 + \sin(3\pi x_{i})) 
      %6 - x_1(1 + \sin(3\pi x_{1})) - x_2(1 + \sin(3\pi x_{2})) 
\end{array}
\right) \ \ \ \ s.t \ \ & x_{1}, x_2 \in [0,1] 
\end{aligned}
\end{equation}
}

%\todo{KT: Ohne multistart probieren, Noch verallgemeinern fuer beliebig viele Dimensionen? Waere interessant, da nichtkonvex!}

%\begin{figure}[H]
%\begin{center}
%\includegraphics[width=0.5\linewidth]{Deb2_paretoFront.png}
%\end{center}
%\caption{Pareto front of test problem \eqref{testproblem_deb2}}
%\label{fig_deb2}
%\end{figure}

\subsection{Implementational details}

\noindent We investigate the reduction in calculation time when using the 
improvements from Section~\ref{sec:improve} by comparing it to a naive implementation from \cite{Tei14}, in which redundant bounds are filtered out afterwards, 
bounds are stored as lists and the maximal box is found by iterating over the nested lists. 
To do the comparison in a fair manner we want to ensure that for all tested instances exactly the same representation points are calculated for both implementations. 
This is not automatically the case because there are typically a lot of boxes of the same size (i.e. the same minimal side length), especially in the beginning stages of the algorithm. For an example consider Figure~\ref{fig_hb_2d}(a) in which the two white boxes have the same minimal side length.
The two implementations might then pick a different box from these equally large boxes depending on the order the boxes are stored internally. 
Therefore, we add the following criteria in case the choice of the box is not unique: 
Firstly, when two boxes have the same minimal side length the box with the greater volume is chosen. %rather than the box with the smaller volume. 
This rules out certain boxes but still not all of them, e.g.\ 
in the example in Figure~\ref{fig_hb_2d}(a) the two white boxes also have the same volume.
Secondly, a kind of ID in form of a specific number is calculated based on the coordinates 
of the lower and upper corner of the box, and the box with the larger number is chosen. This leads to a deterministic choice of the boxes independent of the implementation and, thus, ensures a fair comparison.

%\subsection{Performance improvement when storing the opposing points}
\subsection{Discussion of results for academic test cases}

The results for the three dimensional test problems are shown in Table~\ref{table:N_vs_F_3D}. 
For each test case we indicate a desired approximation quality~$\epsilon$, in our case measured as the maximal minimal side length of all boxes. 
We choose this value relative to the starting box 
which is, in our experiments, defined by the ideal point as lower bound and the nadir point as upper bound. 
The third column indicates the achieved cardinality, i.e. the number of generated nondominated points. 
The next two columns contain the computational time of the implementation from \cite{Tei14} and the improvements discussed in Section~\ref{sec:improve}.
The last column indicates the relative gain in computational time. 
For the highest approximation quality, which is $\epsilon = 0.02$ in our experiments, we reduce the computational time by more than $80\%$ in all five test cases.

\begin{figure}
\centering
\subfloat[Sphere: $\left|\imrepset \right| = 82$]{
  \includegraphics[width=0.32\textwidth]{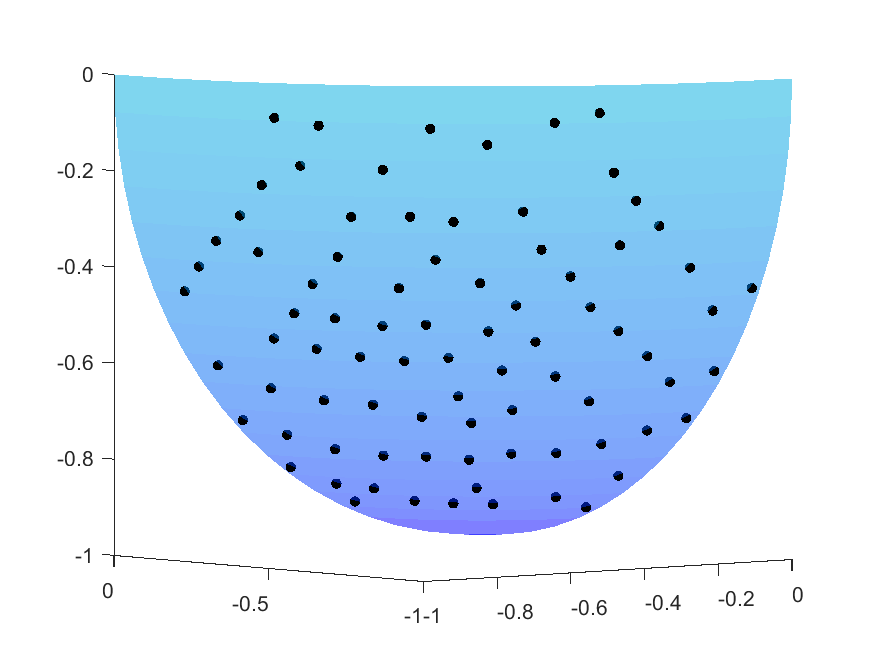}
}
\subfloat[$\left|\imrepset \right| = 316$]{
  \includegraphics[width=0.32\textwidth]{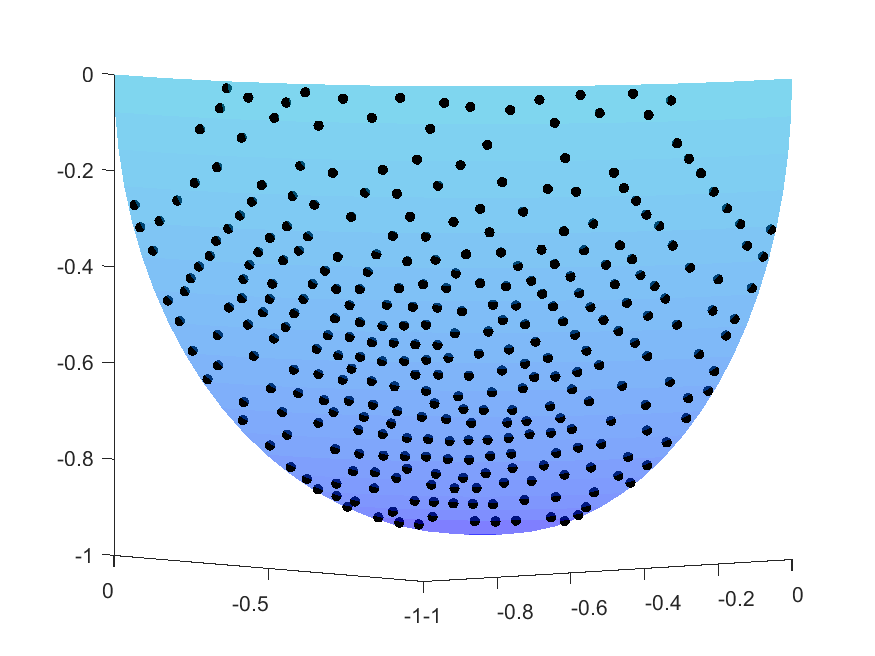}
}
\subfloat[$\left|\imrepset \right| = 1786$]{
  \includegraphics[width=0.32\textwidth]{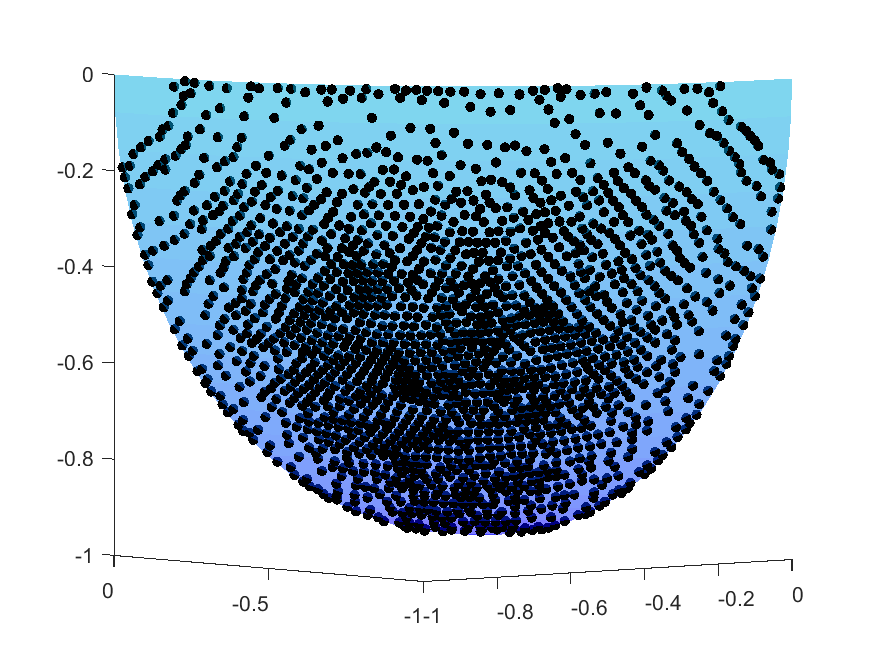}
}
\hspace{0mm}
\subfloat[Ellipsoid: $\left|\imrepset \right| = 132$]{
  \includegraphics[width=0.32\textwidth]{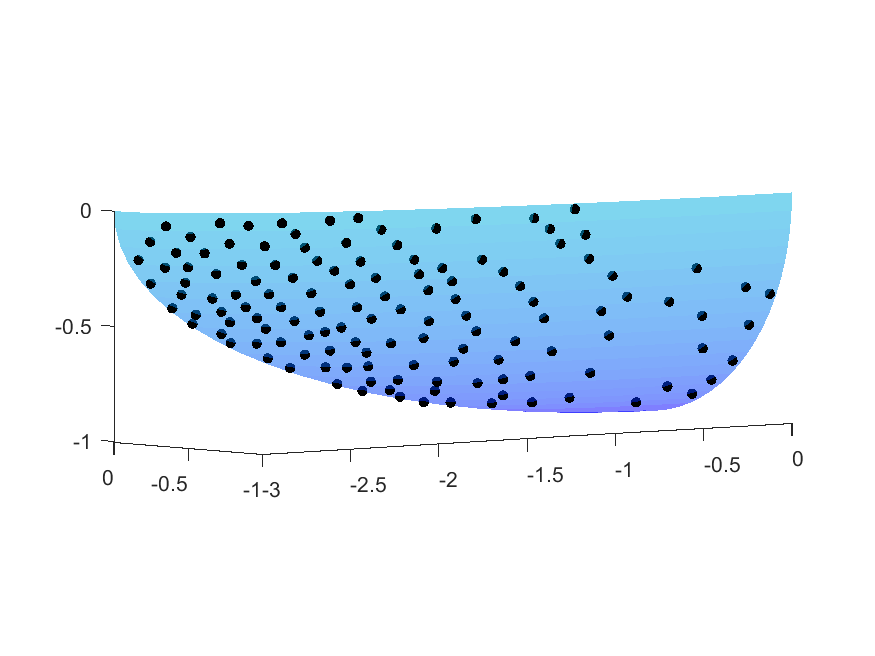}
}
\subfloat[$\left|\imrepset \right| = 498$]{
  \includegraphics[width=0.32\textwidth]{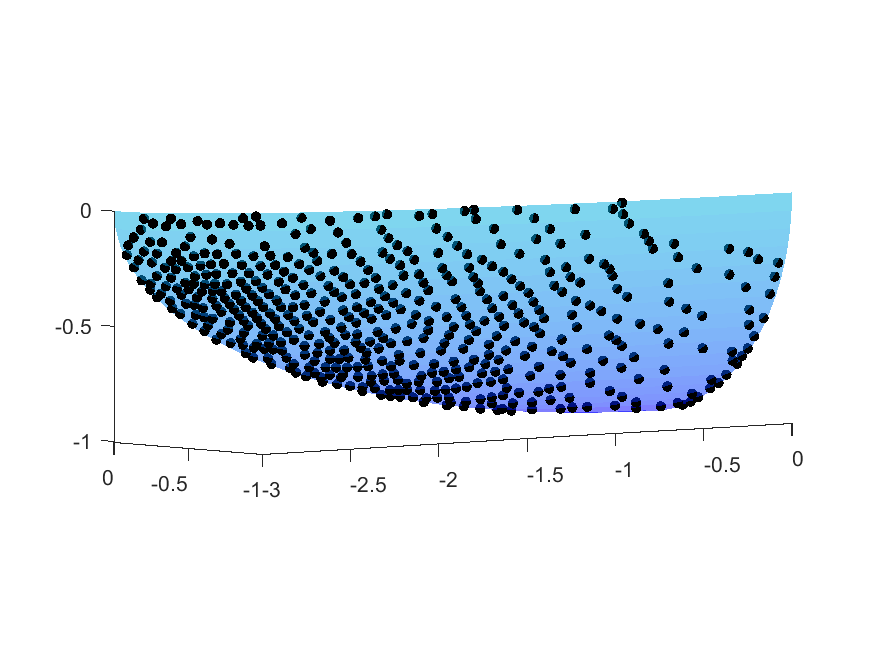}
}
\subfloat[$\left|\imrepset \right| = 2959$]{
  \includegraphics[width=0.32\textwidth]{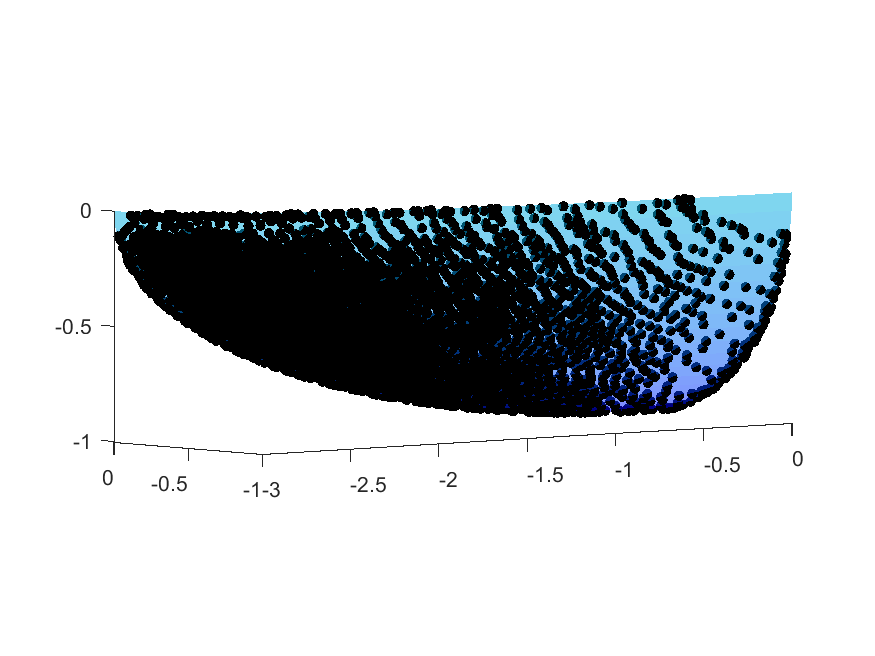}
}
\hspace{0mm}
\subfloat[Non-convex: $\left|\imrepset \right| = 92$]{
  \includegraphics[width=0.32\textwidth]{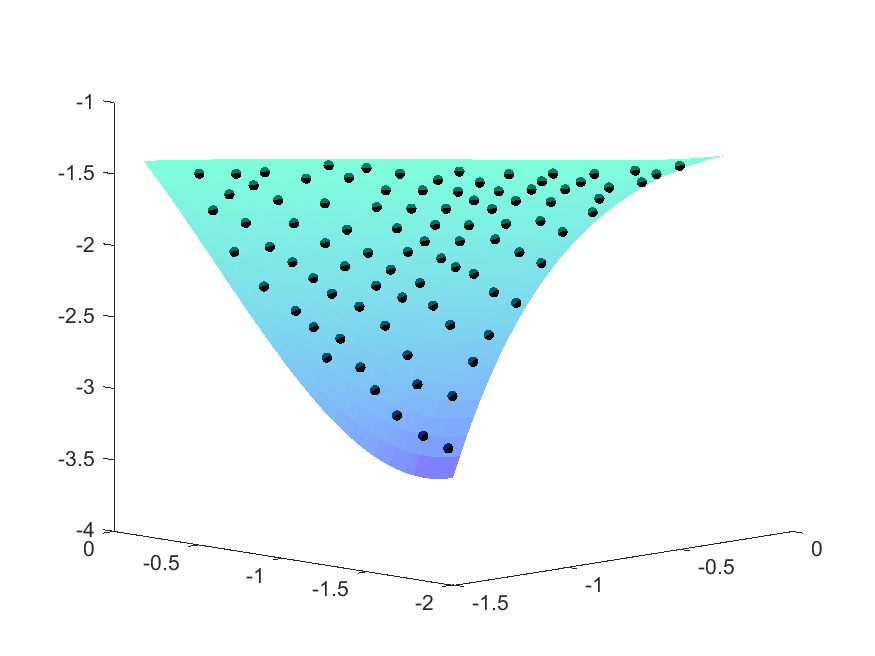}
}
\subfloat[$\left|\imrepset \right| = 380$]{
  \includegraphics[width=0.32\textwidth]{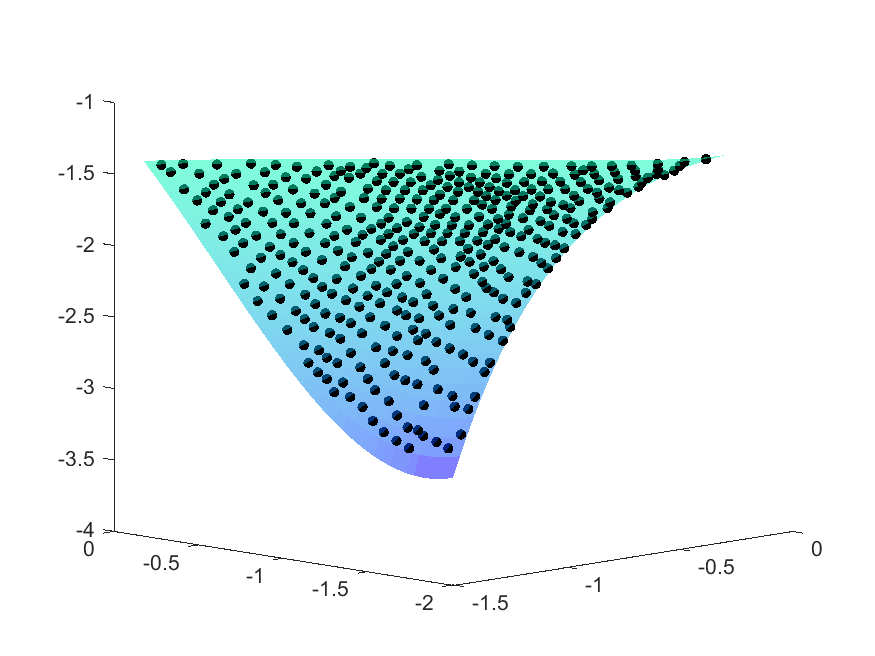}
}
\subfloat[$\left|\imrepset \right| = 2340$]{
  \includegraphics[width=0.32\textwidth]{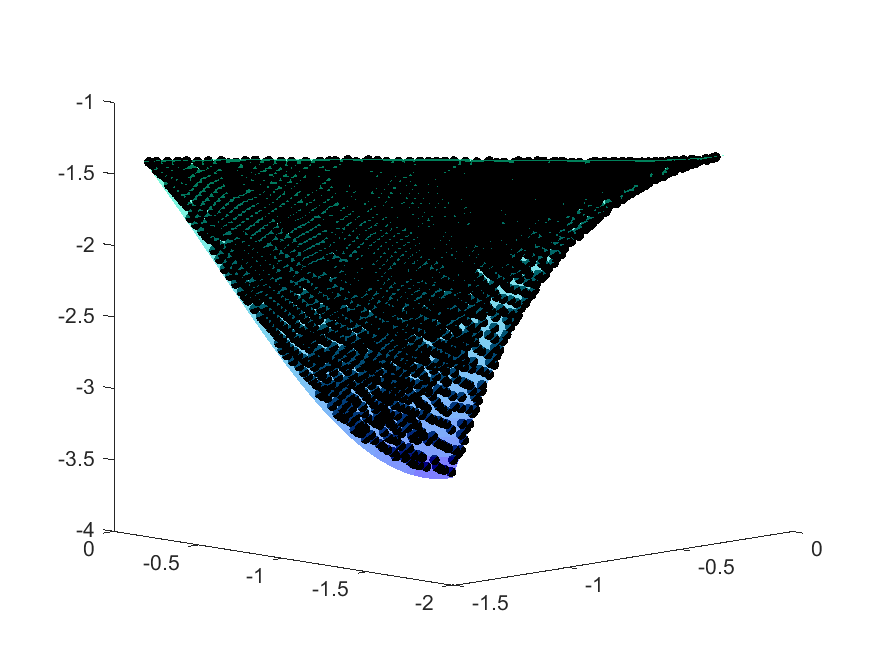}
}
\hspace{0mm}
\subfloat[Comet: $\left|\imrepset \right| = 72$]{
  \includegraphics[width=0.32\textwidth]{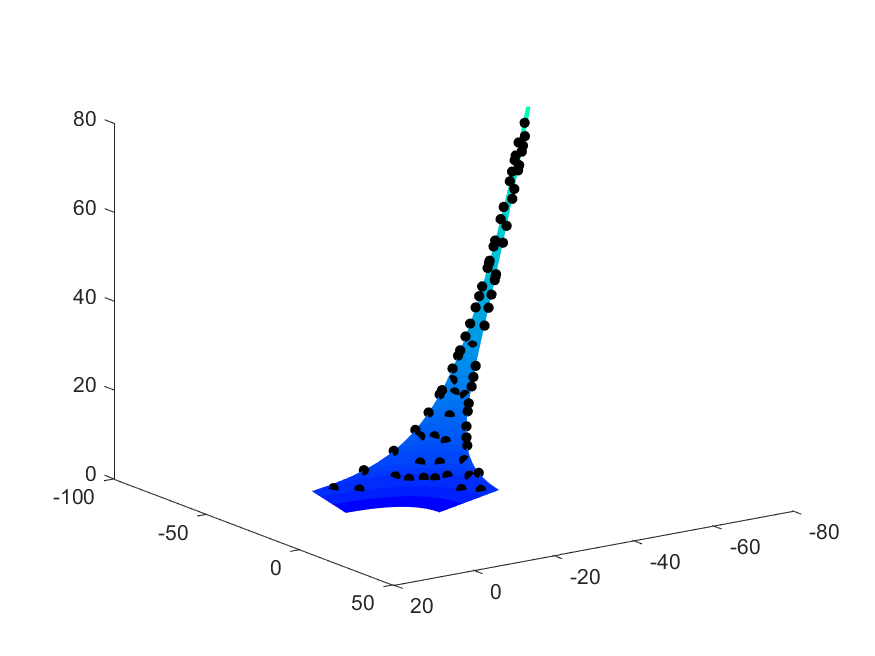}
}
\subfloat[$\left|\imrepset \right| = 326$]{
  \includegraphics[width=0.32\textwidth]{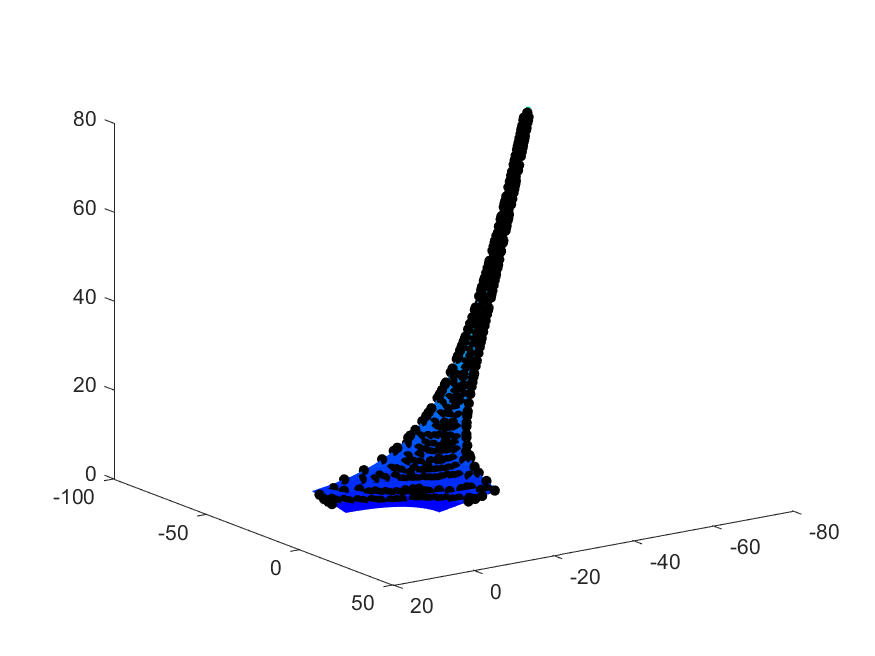}
}
\subfloat[$\left|\imrepset \right| = 946$]{
  \includegraphics[width=0.32\textwidth]{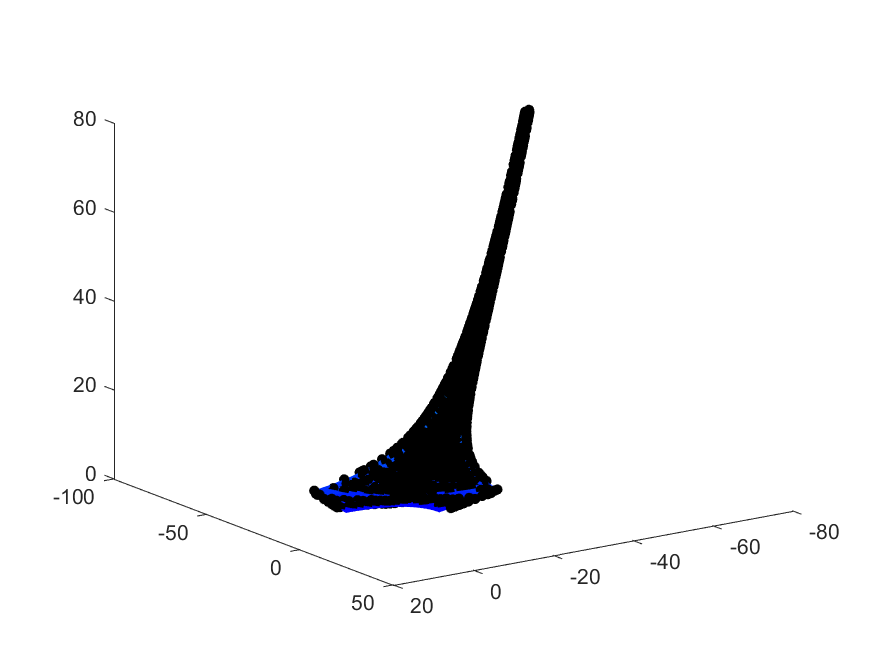}
}
\hspace{0mm}
\subfloat[Patched: $\left|\imrepset \right| = 72$]{
  \includegraphics[width=0.32\textwidth]{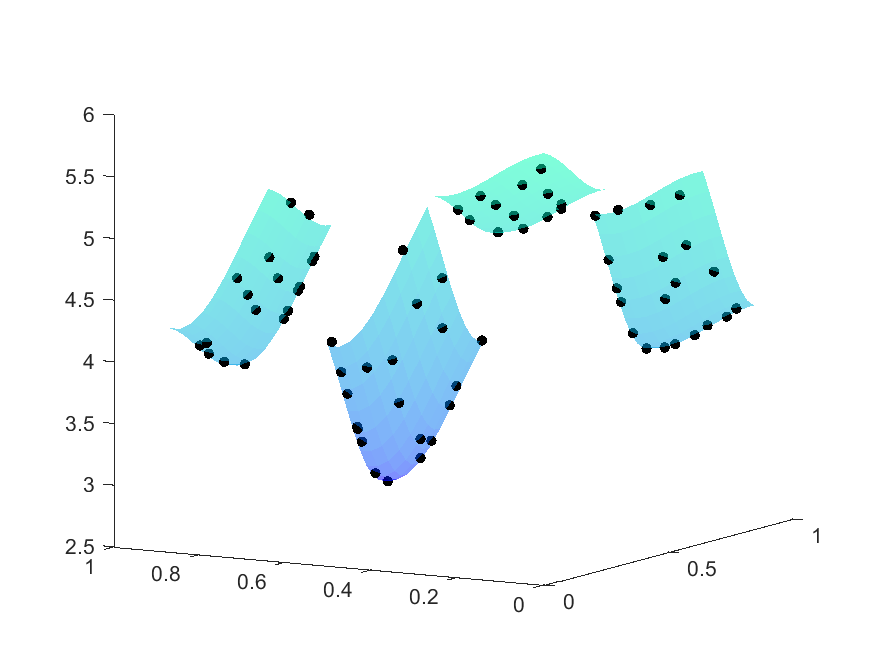}
}
\subfloat[$\left|\imrepset \right| = 237$]{
  \includegraphics[width=0.32\textwidth]{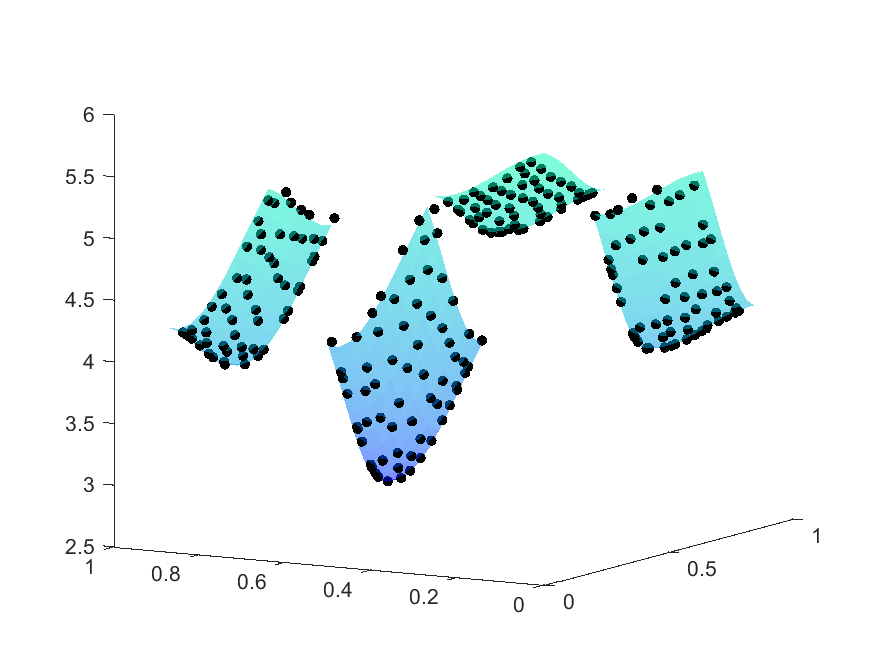}
}
\subfloat[$\left|\imrepset \right| = 1207$]{
  \includegraphics[width=0.32\textwidth]{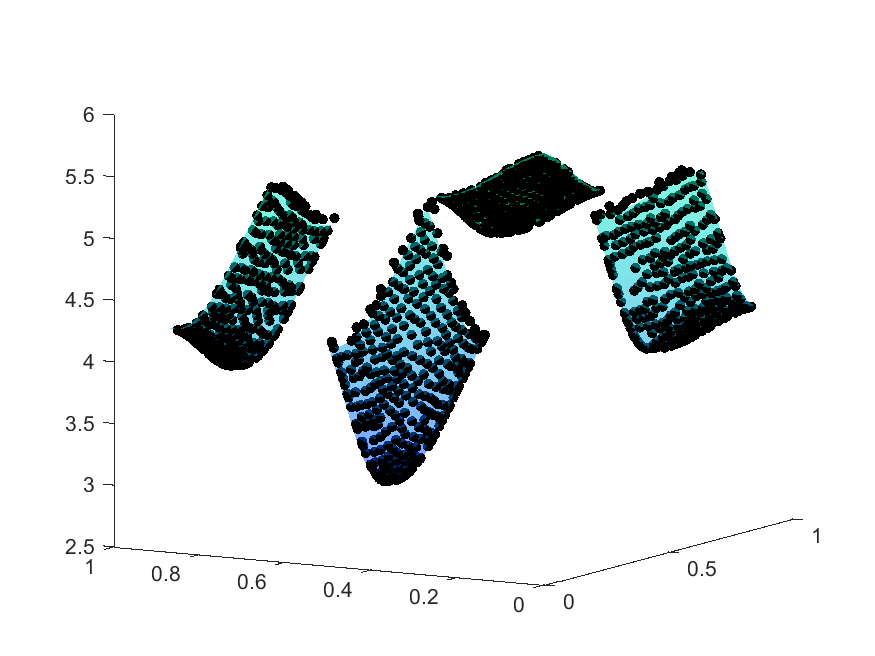}
}
\caption{The representation points for the three dimensional examples as listed in Table~\ref{table:N_vs_F_3D}.}
\label{fig_points_3D}
\end{figure}

The distribution of the resulting representation points for all three-dimensional test problems is displayed in Figure \ref{fig_points_3D}. The higher the approximation quality is, the more accurately the computed representations cover the true Pareto fronts. For some test cases parts at the boundary are represented worse than parts in the center. This is due to the fact that the boxes (that are not drawn in the plots) are rather long and thin and therefore have a small side length. Translated to the plot, the shape of the front is rather straight in these parts and can be approximated well by the generated points. 
%\todo{Hm bzw hat das auch was damit zu tun, ob die Punkte am Rand der Startbox liegen, oder?}

Tables \ref{table:N_vs_F_HS} and \ref{table:N_vs_F_HE} show the results for the hypersphere and hyperellipsoid respectively, over varying dimensions. 
Since the number of boxes does not grow linear with the number of objectives, the same approximation quality implies an increasing computational effort the more objectives we consider. 
As can be seen from the tables, the improved version is vastly superior in terms of calculation time. Even for a non-flat nine-dimensional Pareto front such as the hypersphere, it can produce more than $100$ representation points in 
%reasonable time
less than $500$ seconds, which is impossible with the naive implementation.  

\begin{table}[!htbp]
\centering
\begin{tabular}{l | l | l | c | c | c }

example & $\epsilon$  & 
$\left|\imrepset \right|$
%$\mathcal{N}$  
& $t^{naive}(s)$ & $t^{improved}(s)$ & $\frac{t^{improved}}{t^{naive}}$\\
\hline
\multirow{3}{*}{Sphere} & 0.1  &  82 & 0.016 & 0.001 & 0.06 \\
	                    & 0.05  & 316 & 0.68 & 0.18 & 0.26 \\ 
	                    & 0.02  & 1786 & 87.48  & 8.33 & 0.10 \\
\hline
\multirow{3}{*}{Ellipsoid} & 0.1   & 132 & 0.064  & 0.001 & 0.02 \\
	                    & 0.05   & 498  & 2.28  & 0.44 & 0.19 \\
	                    & 0.02   & 2959 & 383.25  & 27.23 & 0.07 \\
\hline
\multirow{3}{*}{Non-convex} & 0.1    & 92  & 0.032  & 0.001 & 0.03 \\
	                    & 0.05   & 380  & 1.23  & 0.37 & 0.30 \\
	                    & 0.02   & 2340  & 193.07 & 16.55 & 0.09 \\
\hline	                           
\multirow{3}{*}{Comet} & 0.1    & 72 & 0.006  & 0.001 & 0.17 \\
	                    & 0.05   & 326 & 0.75  & 0.25 & 0.33 \\
	                    & 0.02  & 946 & 13.32 & 2.58 & 0.19 \\
\hline	                      
\multirow{3}{*}{Patched} & 0.1    & 72 & 0.018  & 0.001 & 0.06 \\
	                    & 0.05  & 237 & 0.38  & 0.14 & 0.37 \\
	                    & 0.02   & 1207  & 29.15 & 4.33 & 0.15 \\
\end{tabular}
\caption{Comparison of naive and improved implementation for 3D examples}
\label{table:N_vs_F_3D}
\end{table}

\begin{table}[!htbp]
\centering
\begin{tabular}{ l | l | l | c | c | c }

dimension & $\epsilon$ & 
%$\left|\imrepset \right|$  
$\left|\imrepset \right|$
& $t^{naive}(s)$ & $t^{improved}(s)$ & $\frac{t^{improved}}{t^{naive}}$ \\
\hline
\multirow{3}{*}{4} &  0.2   &  87   &  0.237  &  0.036 & 0.15 \\
	               &  0.15  &  181  &  1.6    &  0.19  & 0.12 \\
	               &  0.1   &  638  &  56.45  &  3.78  & 0.07 \\
\hline
\multirow{3}{*}{5} & 0.3  &  51   &  0.40  &  0.04   & 0.10 \\
	               & 0.2  &  271  &  94.35 &  2.75   & 0.03 \\
	               & 0.15 &  816  &  $--$  &  36.206 & $--$ \\
\hline	                           
\multirow{3}{*}{6} & 0.35  &  41   &  1.41  &  0.08  & 0.06 \\
	               & 0.3   &  96   &  48.78 &  0.85  & 0.02 \\
	               & 0.25  &  260  &  $--$  &  11.43 & $--$ \\
\hline	                      
\multirow{3}{*}{7} & 0.35  &  45   &  10.92  &  0.37  & 0.03 \\
	               & 0.325 &  84   &  327.82 &  2.17  & 0.007 \\
	               & 0.3   &  159  &  $--$   &  13.52 & $--$  \\
\hline	                      
\multirow{2}{*}{8} & 0.35  &  56   &  120.45 &  1.41  & 0.01 \\
	               & 0.325 &  152  &  $--$   &  50.32 & $--$ \\
\hline	                      
\multirow{1}{*}{9} & 0.325 &  157  &  $--$  &  179.62 & $--$ \\
\hline	                      
\end{tabular}
\caption{Comparison of naive and improved implementation for the hypersphere. Solution times higher than $500$ seconds are marked by $--$.}
\label{table:N_vs_F_HS}
\end{table}

\begin{table}[!htbp]
\centering
\begin{tabular}{ l | l | l | c | c | c }

dimension & $\epsilon$ & 
%$\left|\imrepset \right|$  
$\left|\imrepset \right|$
& $t^{naive}(s)$ & $t^{improved}(s)$ & $\frac{t^{improved}}{t^{naive}}$ \\
\hline
\multirow{3}{*}{4} &  0.2   &  157   &  1.1    &  0.14  &  0.13 \\
	               &  0.15  &  363   &  10.66  &  0.98  &  0.09 \\
	               &  0.1   &  1197  &  323.21 &  13.95 &  0.04 \\
\hline
\multirow{3}{*}{5} & 0.3  &  101   &   4.13  &  0.2   & 0.05 \\
	               & 0.2  &  552   &   $--$  &  12.73 & $--$ \\
	               & 0.15 &  1747  &   $--$  &  198.2 & $--$ \\
\hline	                           
\multirow{3}{*}{6} & 0.35  &  74    &   11.09  &  0.30 & 0.03 \\
	               & 0.3   &  167   &   254.53 &  3.01 & 0.01 \\
	               & 0.25  &  517   &   $--$   &  55.9 & $--$ \\
\hline	                      
\multirow{3}{*}{7} & 0.35  &  116   &   $--$  &  4.92   & $--$ \\
	               & 0.325 &  209   &   $--$  &  24.79  & $--$ \\
	               & 0.3   &  369   &   $--$  &  107.88 & $--$ \\
\hline	                      
\multirow{2}{*}{8} & 0.35  &  96   &  $--$  &  5.99   & $--$\\
	               & 0.325 &  274  &  $--$  &  218.66 & $--$\\
\hline	                      
\multirow{1}{*}{9} & 0.325 &  215  &  $--$  &  252.6 & $--$ \\
\hline	                      
\end{tabular}
\caption{Comparison of naive and improved implementation for the hyperellipsoid. Solution times higher than $500$ seconds are marked by $--$.}
\label{table:N_vs_F_HE}
\end{table}

\subsection{Real world application: radiotherapy planning}

\begin{figure}[ht]
\centering
\subfloat[A multi-leaf collimator. IMRT planning seeks to determine the best way of steering the collimator leaves in and out of the radiation field such that the resulting dose is optimal for the given patient geometry.]{
  \includegraphics[width=0.44\textwidth]{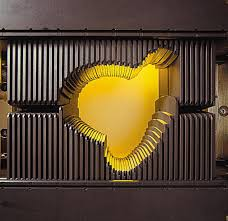}
}
\qquad
\subfloat[The resulting dose distribution. Over a horizontal CT slice of the patient anatomy, regions of high dose are colored red and regions of low dose are colored blue.]{
  \includegraphics[width=0.45\textwidth]{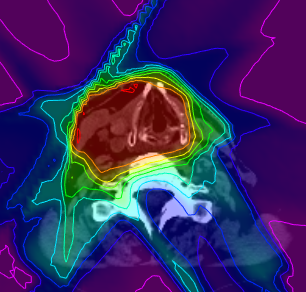}
}
\caption{}
\label{fig_IMRT}
\end{figure}

In intensity-modulated radiation therapy (IMRT), the patient is irradiated with photon beams to destroy the cancer cells. The profile of each photon beam can be modulated by moving the leaves of a multi-leaf collimator (Figure \ref{fig_IMRT}(a)) in and out of the beam, thus partially blocking it. Utilizing these degrees of freedom, a radiation dose distribution can be achieved that conforms very well to the tumor shape (Figure \ref{fig_IMRT}(b)). 

Finding the best beam profiles $x$ in the set of possible profiles $X$ is a multi-objective optimization problem. Some objectives and constraints measure how well the delivered dose $d(x)$ conforms to the tumor, thus predicting the chance for curing the cancer. Other objectives and constraints measure the dose delivered to specific organs or regions at risk near the tumor, thus controlling the risk for specific side effects of the therapy.   

In what is called the a-posteriori approach to solving the IMRT planning problem, the first and computationally expensive step is to find reasonably many solutions that are evenly distributed along the Pareto front of the IMRT planning problem (see \cite{KueHam02}, \cite{BokFor12}). In the context of radiotherapy planning, the resulting representative set is called the plan database. 

\begin{table}[b]
\centering
\begin{tabular}{| l | l |}
\hline
structure & function \\
\hline
\multicolumn{2}{c}{\textbf{objectives}}  \\
\hline
  CTV   &  $f^{\mathrm{CTV}} = \sum_{i \in I^{\mathrm{CTV}} } (d_i(x) - 60)^2$  \\
  Parotid left  &  $f^{\mathrm{Pl}} = \sum_{i \in I^{\mathrm{Pl}} } d_i(x)^2$  \\
  Parotid right &  $f^{\mathrm{Pr}} = \sum_{i \in I^{\mathrm{Pr}} } d_i(x)^2$  \\
  Myelon  & $f^{\mathrm{M}} = \sum_{i \in I^{\mathrm{M}} } d_i(x)^2$  \\
  Eye left & $f^{\mathrm{El}} = \sum_{i \in I^{\mathrm{El}} } d_i(x)^2$   \\
  Eye right & $f^{\mathrm{Er}} = \sum_{i \in I^{\mathrm{Er}} } d_i(x)^2$   \\
\hline
\multicolumn{2}{c}{\textbf{constraints}}  \\
\hline
  CTV   &  $g_1^{\mathrm{CTV}} = \sum_{i \in I^{\mathrm{CTV}} } \max \{0, 55 - d_i(x) \} \leq 0$  \\
  CTV   &  $g_2^{\mathrm{CTV}} = \sum_{i \in I^{\mathrm{CTV}} } \max \{0, d_i(x) - 66 \} \leq 0$   \\
  Myelon &  $g^{\mathrm{M}} = \sum_{i \in I^{\mathrm{M}} } \max \{0, d_i(x) - 45 \} \leq 0$    \\
\hline
\end{tabular}
\caption{Objectives and constraints for the IMRT problem.}
\label{table:IMRT_obj_and_cons}
\end{table}

We applied the hyperboxing algorithm to compute the plan database for a head and neck case. 
The objectives and constraints used for this case are listed in Table \ref{table:IMRT_obj_and_cons}. 

In order to handle the IMRT problem numerically, the patient anatomy is voxelized, and the dose $d(x)$ is given as a vector $(d_{i}(x))_{i \in I}$ with each entry denoting the dose inside a voxel. The objectives and constraints are mathematically modelled as one- or two-sided quadratic deviations of the voxel dose values $d_i$ from a reference value (often equal to zero), summed over all voxels $I^s$ of the corresponding structure $s$. 

For this problem, we choose the box $[0, 3000]^6$ as starting box rather than approximating the nadir. Plans outside this box, while being Pareto optimal, are not interesting for a clinician as one or more objectives are not good enough for the plan to be acceptable.

\begin{figure}[ht]
\centering
\subfloat[]{
  \includegraphics[width=0.33\textwidth]{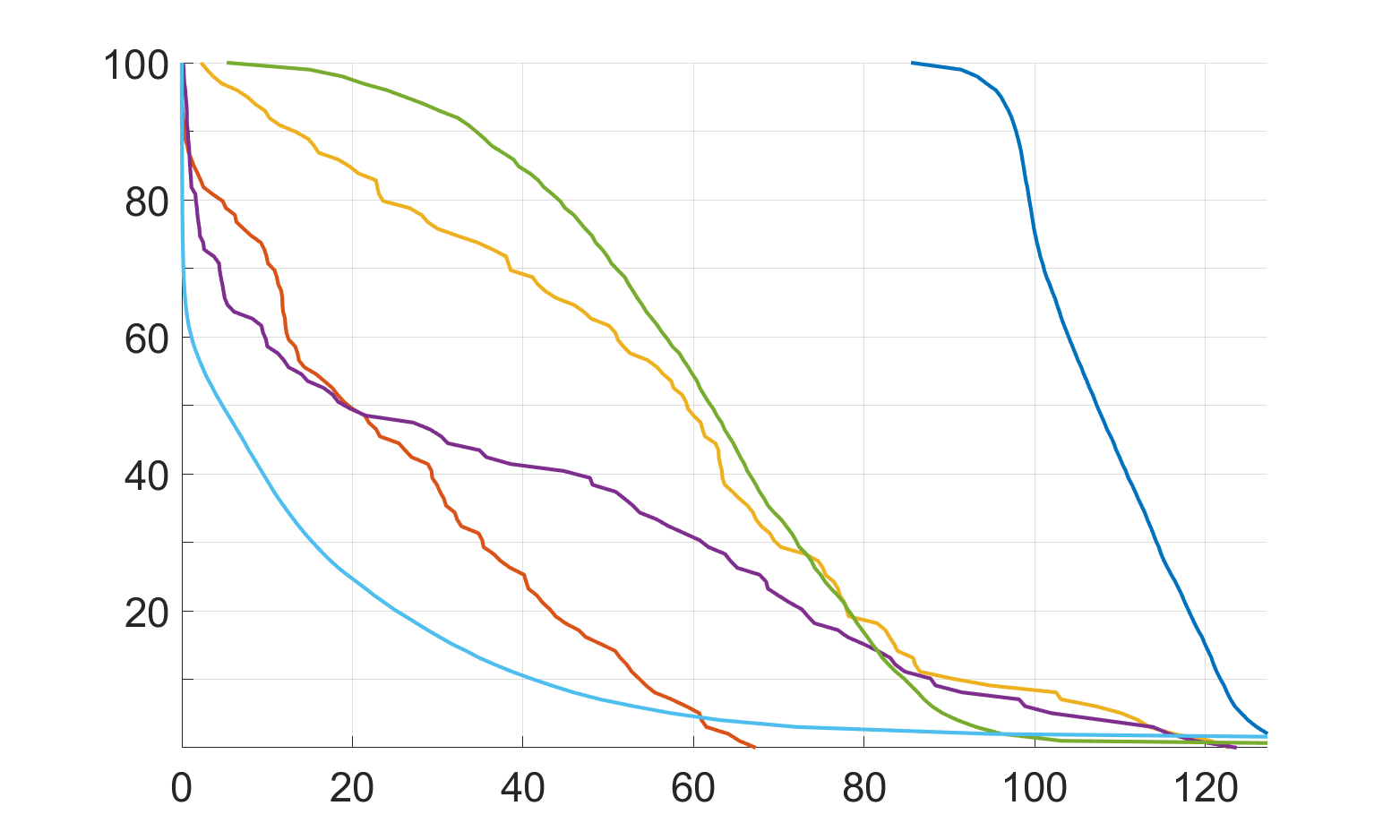}
}
\subfloat[]{
  \includegraphics[width=0.33\textwidth]{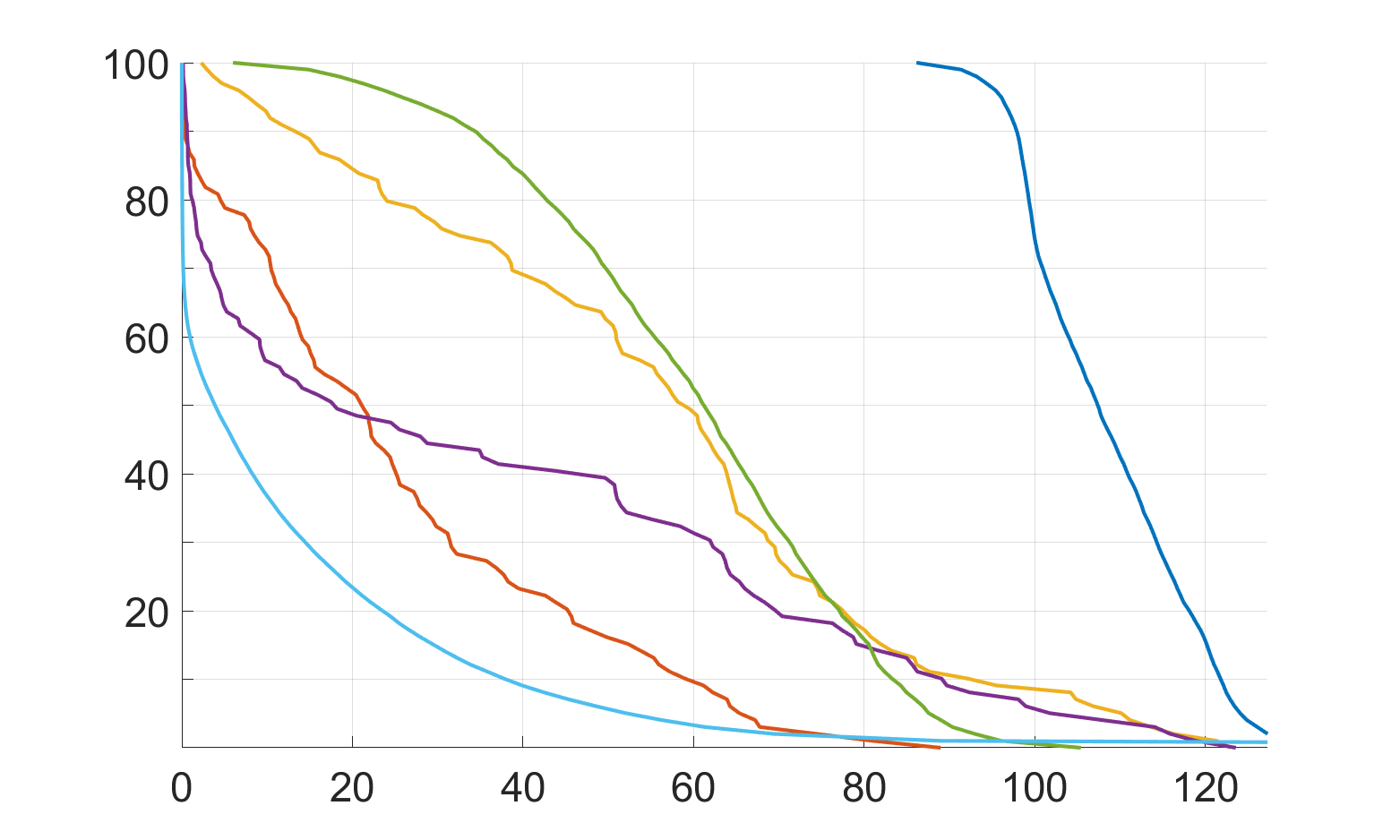}
}
\subfloat[]{
  \includegraphics[width=0.33\textwidth]{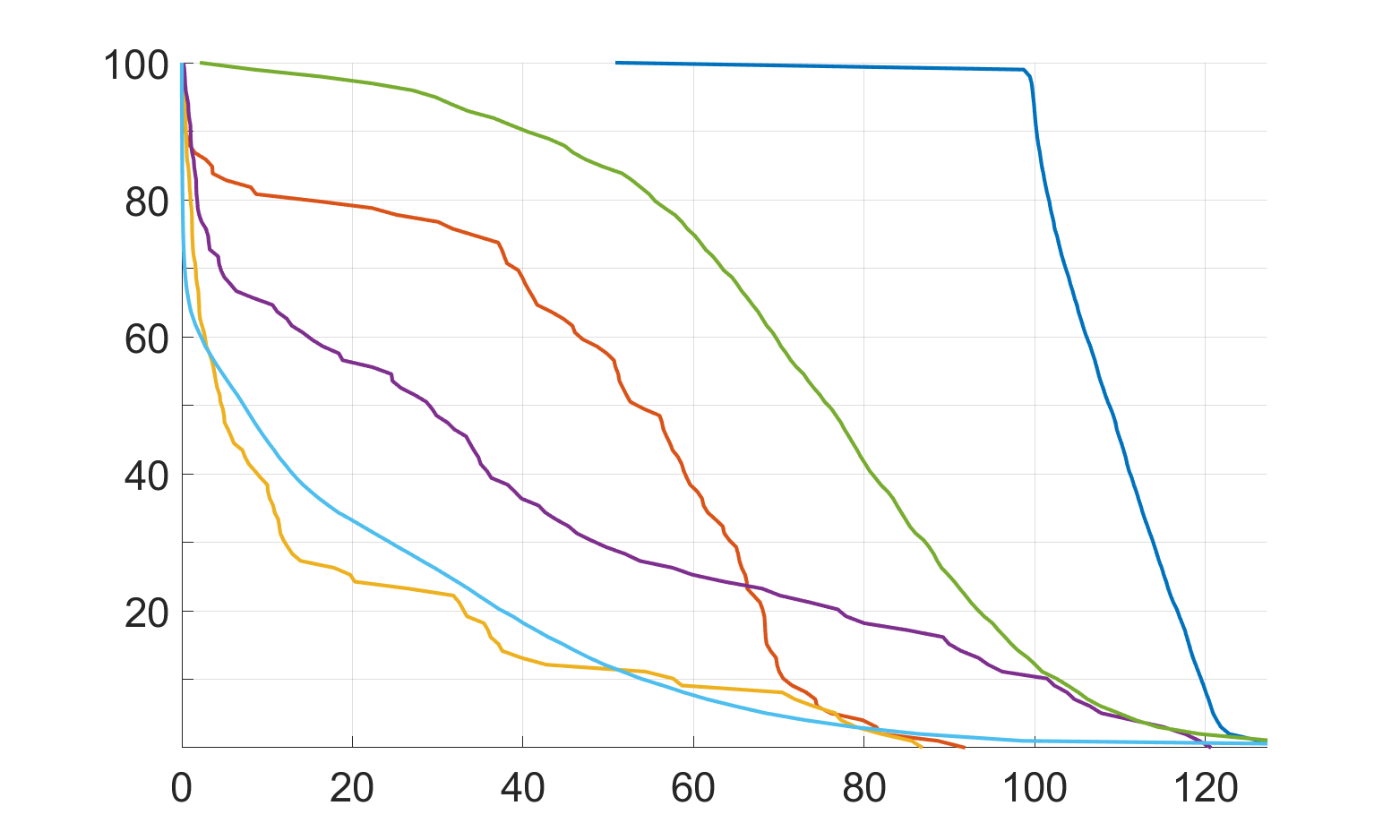}
}
\hspace{0mm}
\subfloat[]{
  \includegraphics[width=0.33\textwidth]{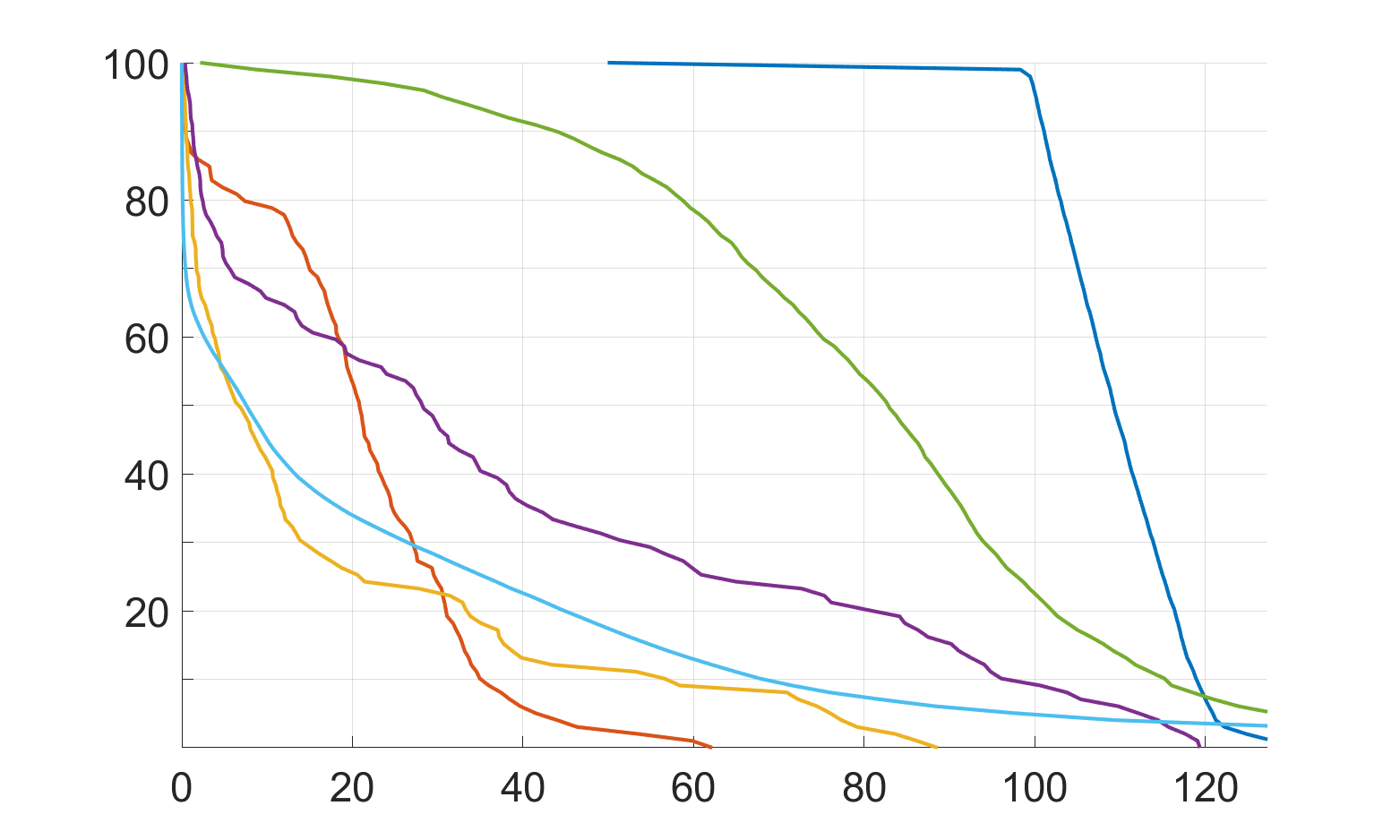}
}
\subfloat[]{
  \includegraphics[width=0.33\textwidth]{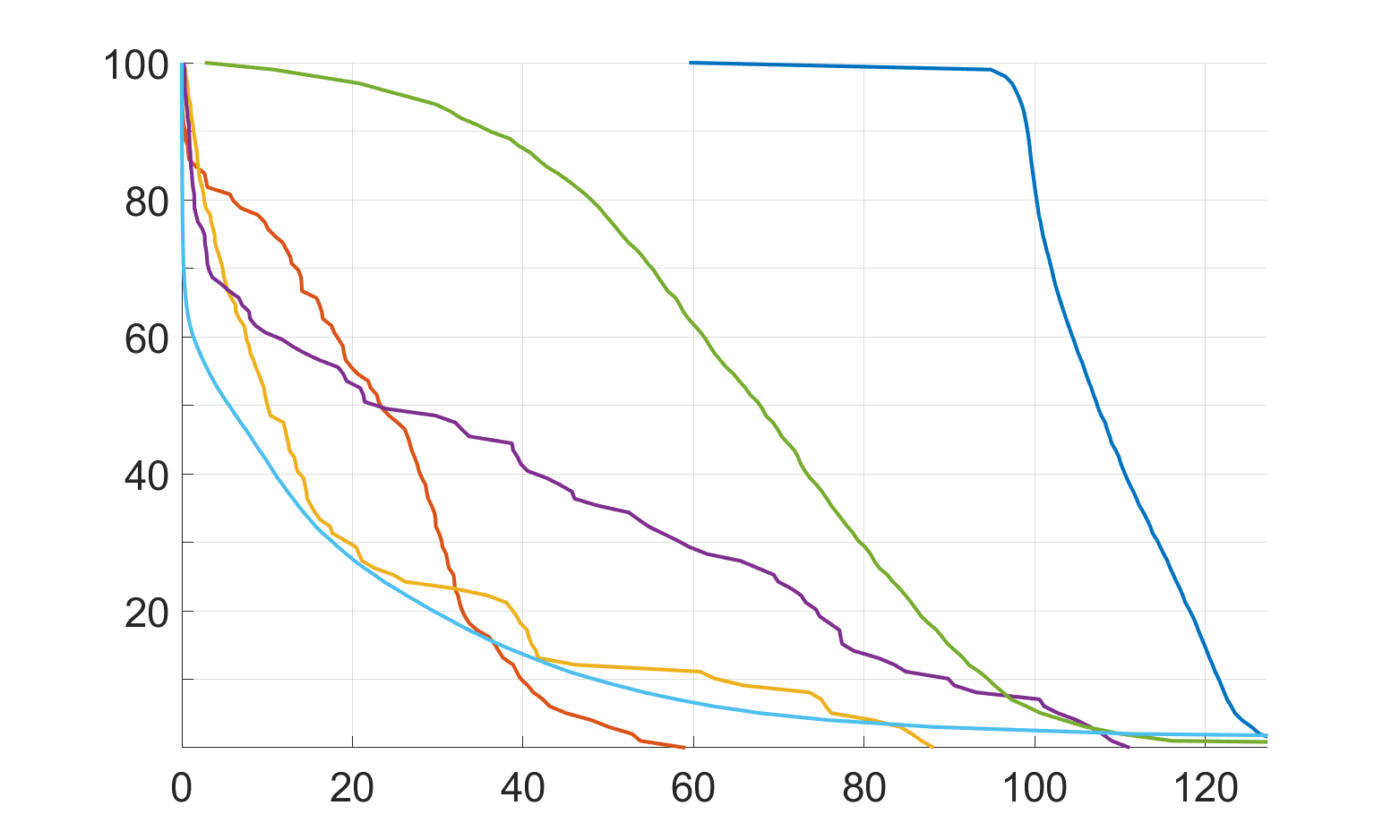}
}
\subfloat[]{
  \includegraphics[width=0.33\textwidth]{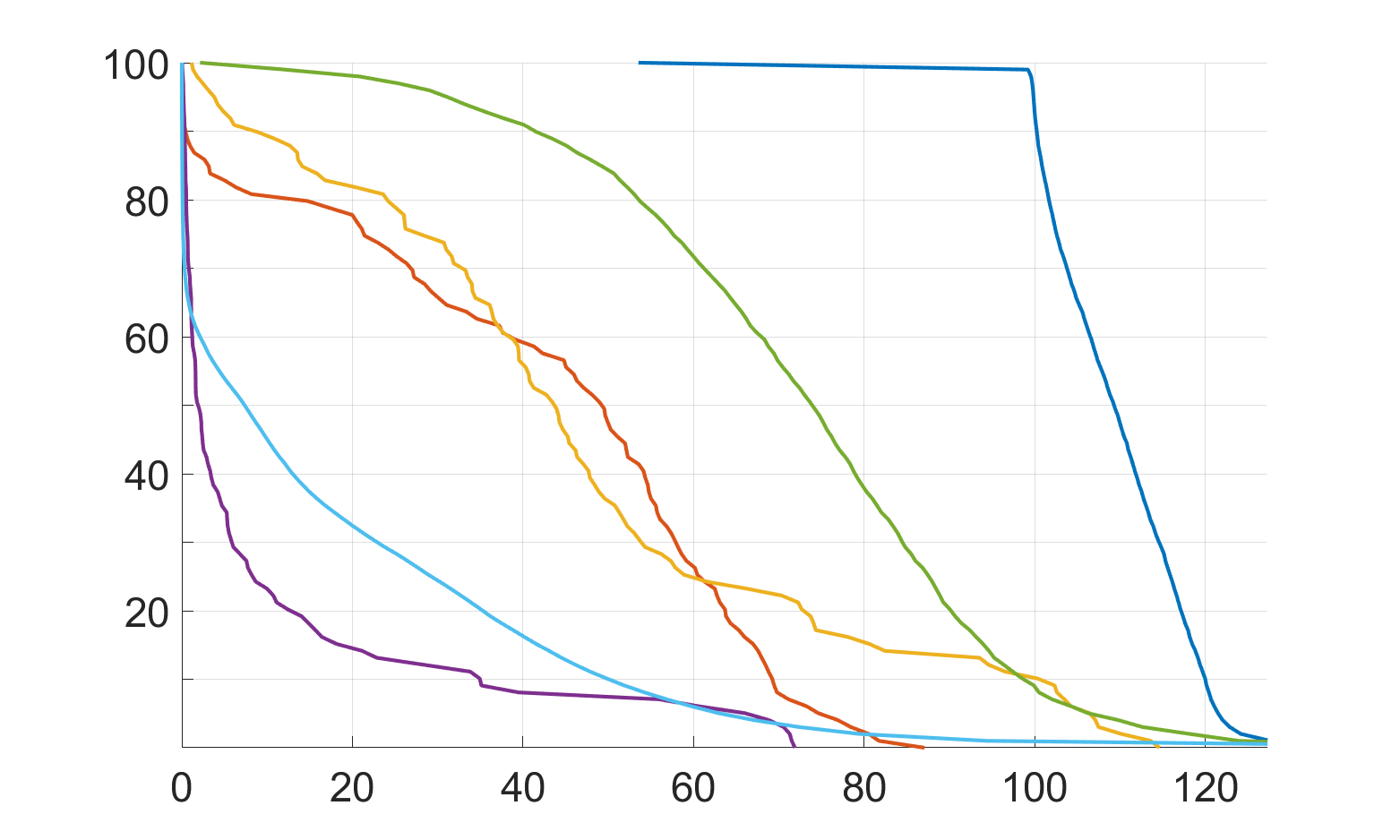}
}
\hspace{0mm}
\subfloat[]{
  \includegraphics[width=0.33\textwidth]{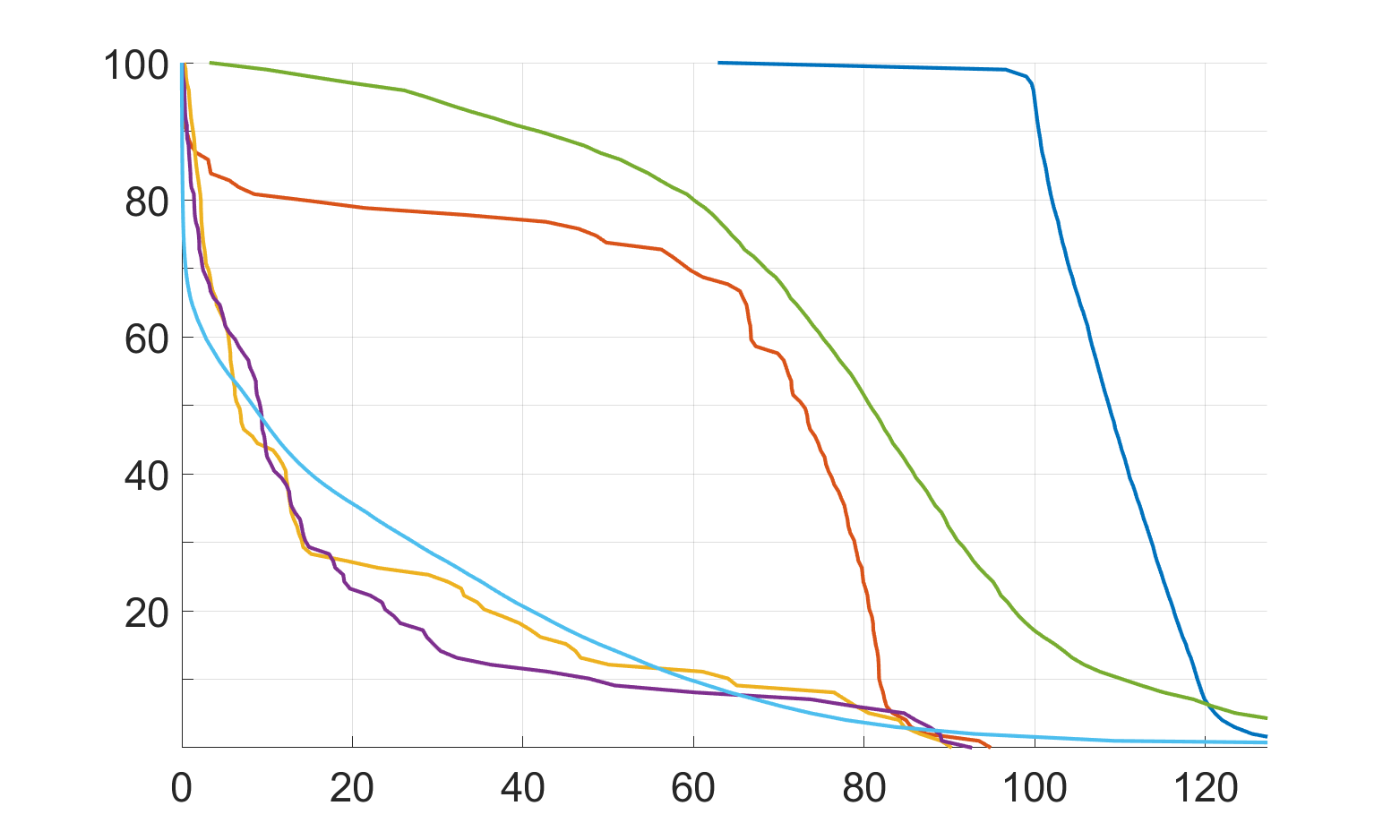}
}
\subfloat[]{
  \includegraphics[width=0.33\textwidth]{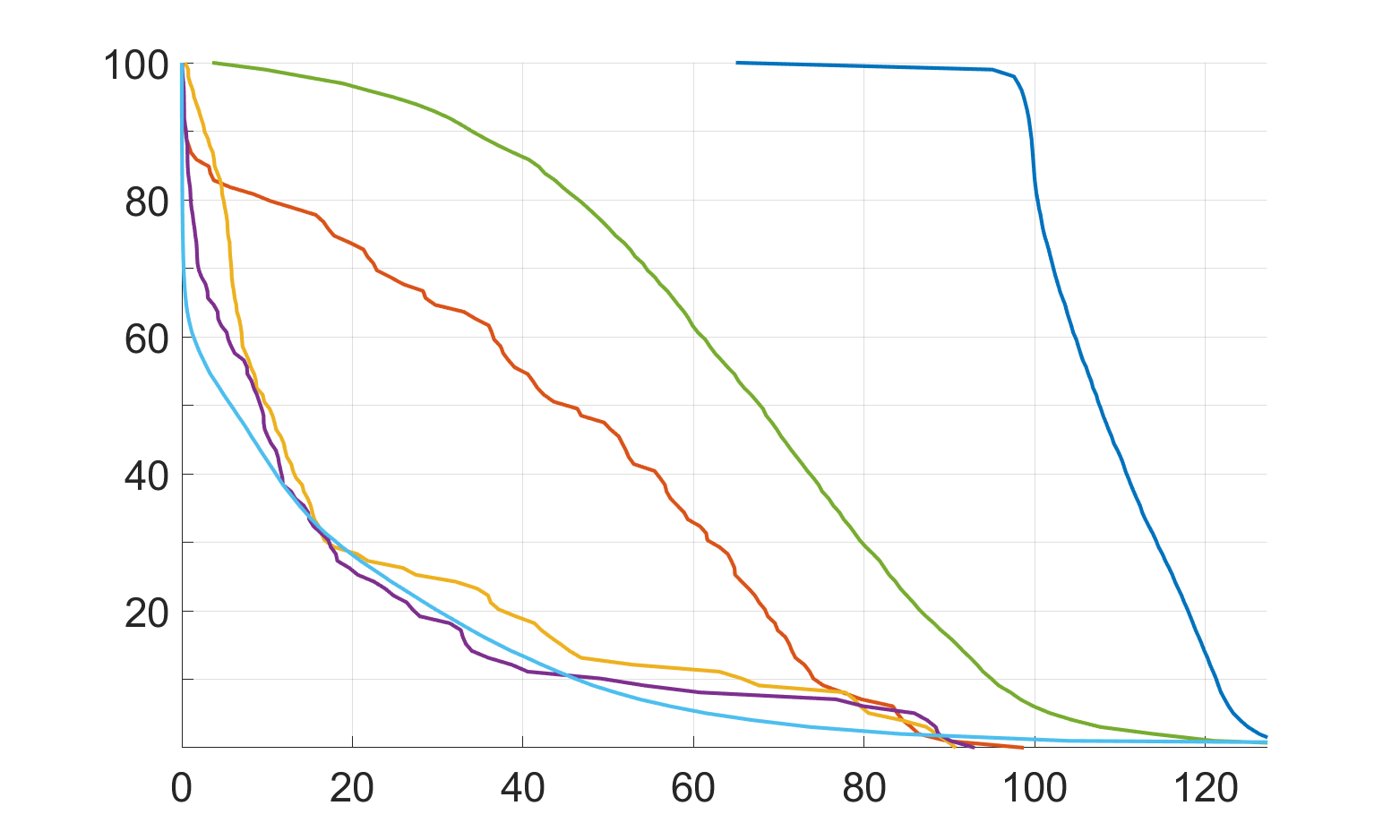}
}
\subfloat[]{
  \includegraphics[width=0.33\textwidth]{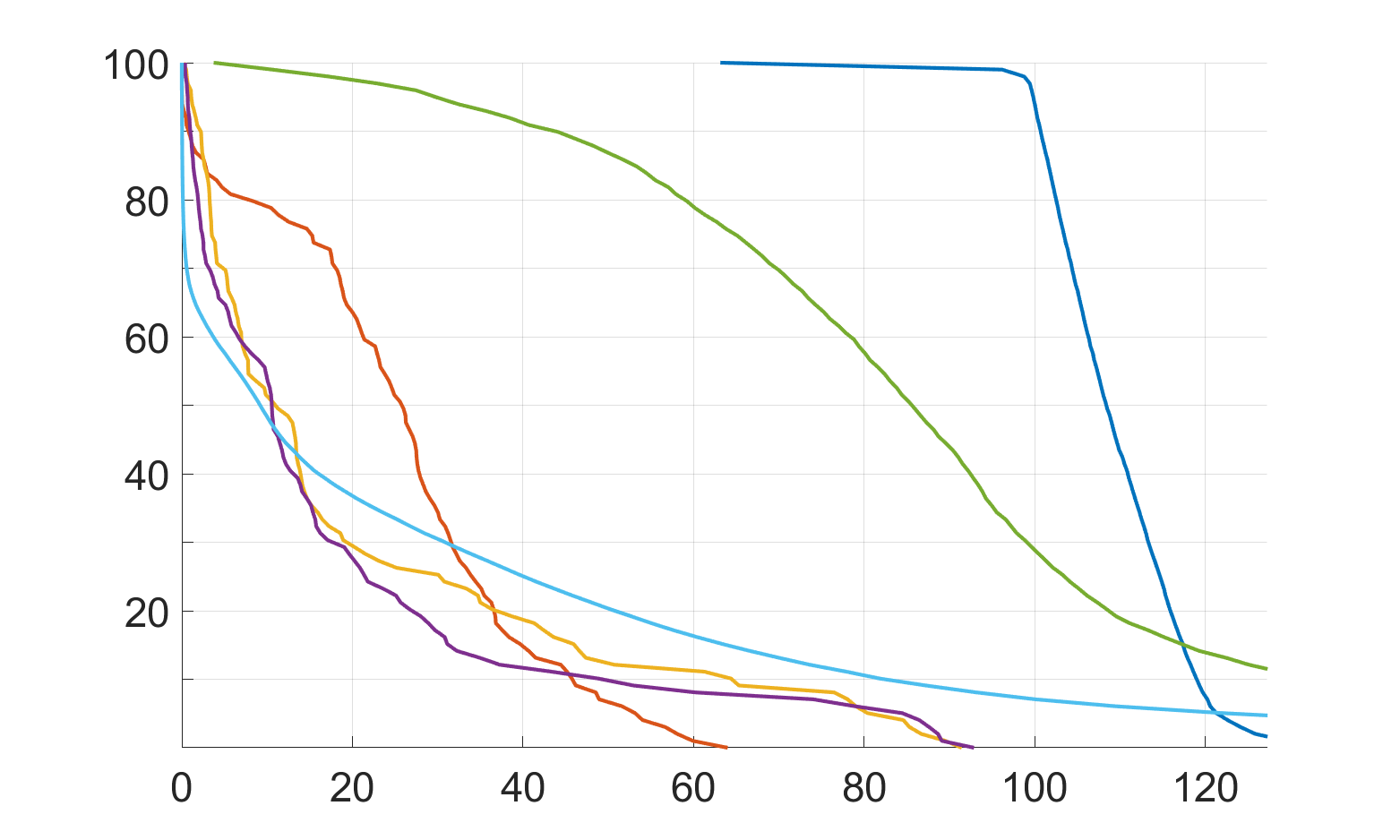}
}
\hspace{0mm}
\subfloat[]{
  \includegraphics[width=0.33\textwidth]{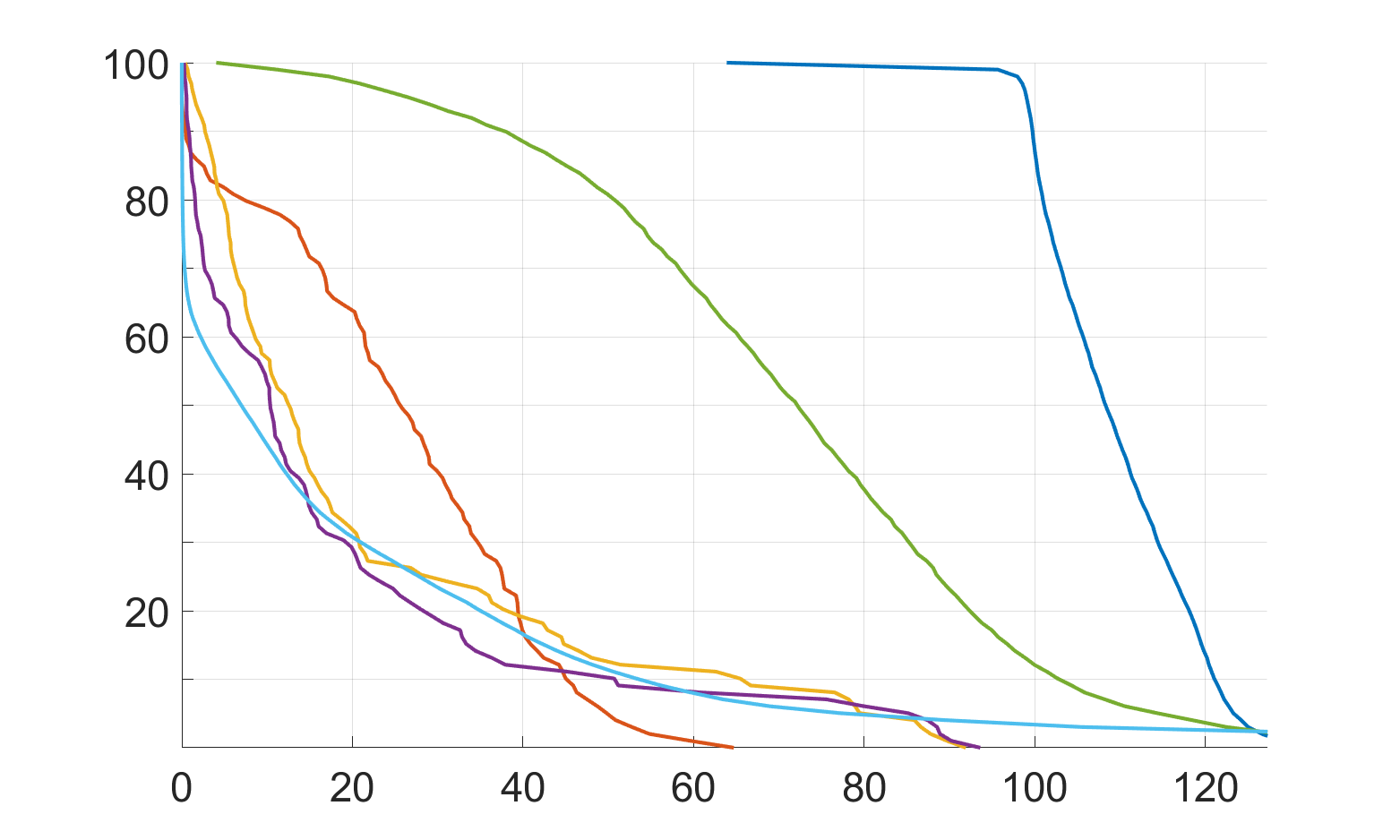}
}
\subfloat[]{
  \includegraphics[width=0.33\textwidth]{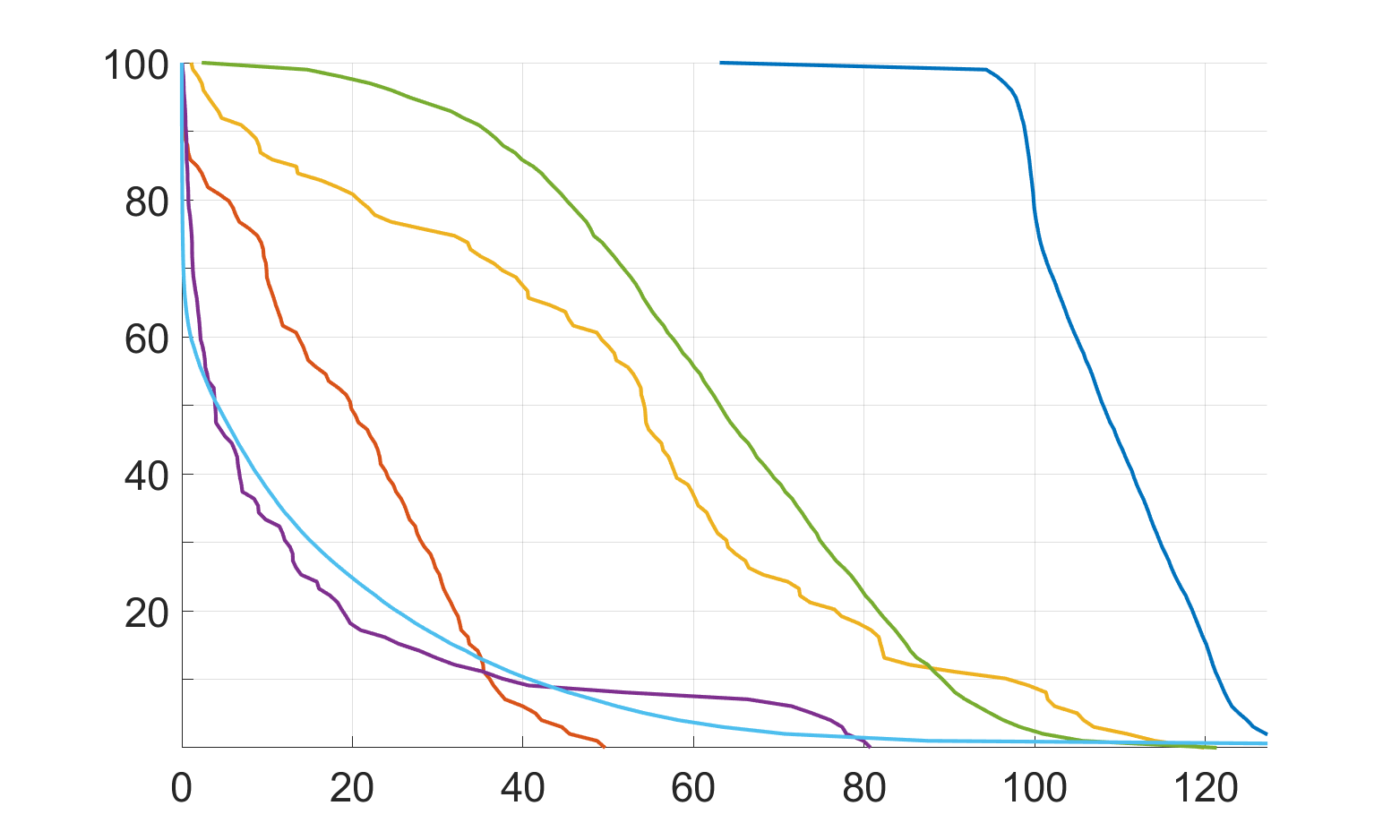}
}
\subfloat[]{
  \includegraphics[width=0.33\textwidth]{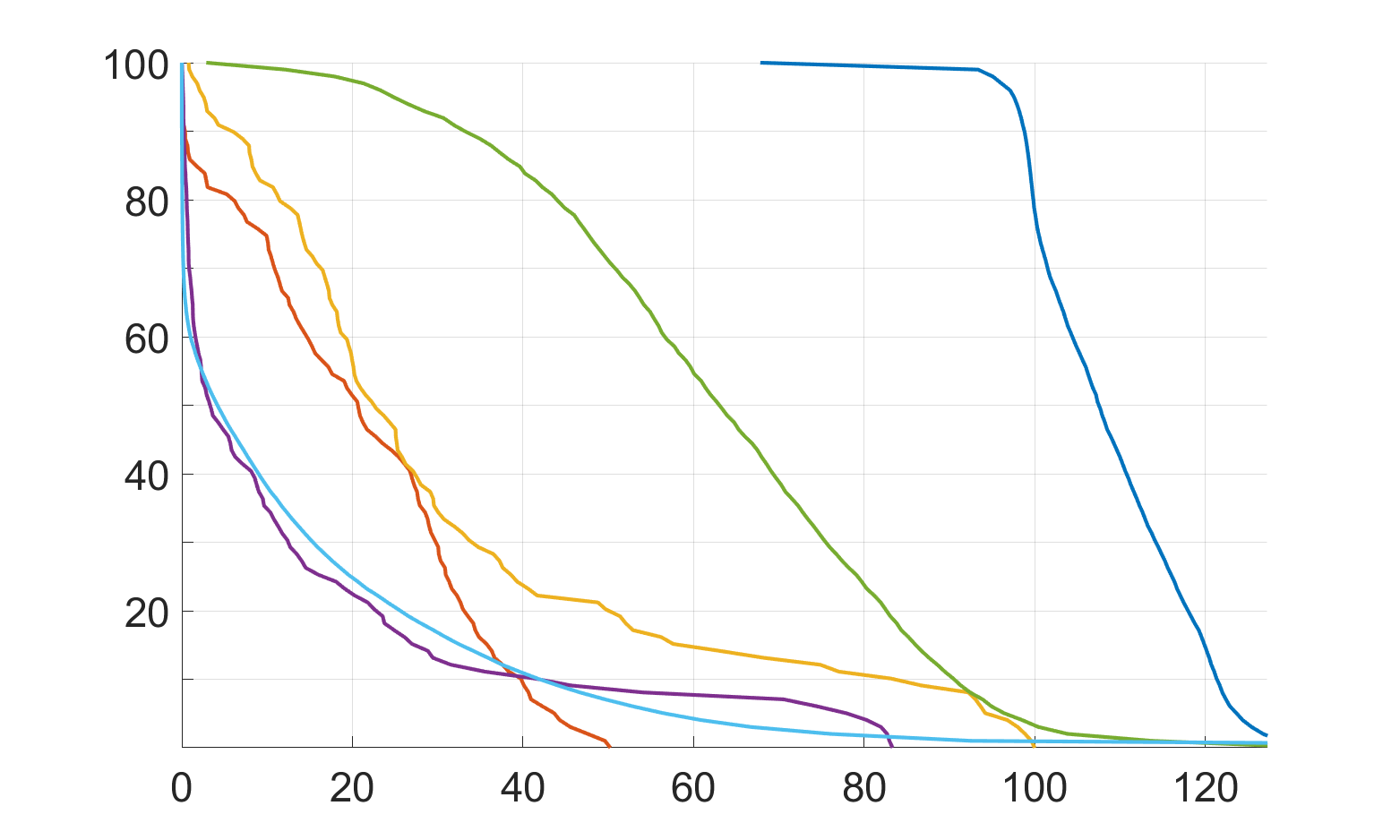}
}
\caption{The DVH diagrams of the first 12 plans for the head neck IMRT planning problem. The structures are CTV (dark blue), right parotid (green), left parotid (purple), left eye (red), right eye (yellow) and myelon (light blue).}
\label{fig_IMRT_dvhs}
\end{figure}

Figure~\ref{fig_IMRT_dvhs} shows the first 12 database plans as DVH (dose volume histogram) diagrams. 
A DVH diagram plots the volume percentage of a structure (y-axis) that receives at least a certain amount of dose (x-axis, here given as percentage of the prescribed dose of 60 Gy). There is one curve for each of the six
%https://www.overleaf.com/project/5c407568bb75b013682072e8 
structures: target volume (CTV), left parotid, right parotid, myelon, left eye and right eye. 
Note that the target volume (dark blue) should receive a high amount of dose while all the other five structures should be preserved as much as possible. 

As can be seen from these DVH diagrams, the hyperboxing algorithm produces systematically varied plans, each with distinctive trade-offs between the six objectives. 
While the goal of hitting the target volume is achieved in all plans, the surrounding structures are affected differently in the presented plans. 
In order to reduce the box size to $0.0625$ of the original size, $144$ plans were calculated. The hyperboxing algorithm took $3.681$ seconds for generating the representation, while the individual optimization runs to solve the Pascoletti-Serafini problems took much longer ($\approx 1$ minute per plan).   

Once the plan database is computed, the database plans can be linearly interpolated to attain a continuous approximation of the Pareto front. The planner can then interactively steer the interpolation coefficients using a slider panel in order to find the best compromise (see \cite{MonKue08}). It has been shown that MCO is capable to both shorten the planning time and improve the plan quality (see \cite{CraHon12}, \cite{KieVis15}).

\section{Conclusion and Outlook} \label{sec:concl}
In this paper we presented an algorithmic improvement of the class of hyperboxing algorithms that makes them applicable to problems with a higher number of objectives. The latter is important for real-world applications which often have more than three objectives. Especially when a decision maker is involved, algorithms with a good performance that respond rather quickly are necessary.
In our tests the proposed algorithmic improvement reduces the computational time considerably
while maintaining the quality of the representation, in our case a certain maximal coverage error. If the decision maker wants to specify a desired cardinality instead, the algorithm at hand can be quickly adapted to stop as soon as a given number of representative points is reached. 

In our numerical tests we chose a maximal coverage error between $0.1$ and $0.35$. However, before knowing anything about the Pareto front of the problem at hand, it is difficult to select this value properly. In some cases, a rather high coverage error would be completely fine while in other cases a finer representation of the Pareto front is required. Since the computational time increases drastically, especially with many objective functions, one further improvement could be to reduce the desired maximal coverage error in a step-wise manner. This means to start with a rather high value, e.g., $\epsilon = 0.5$. 
If the representation is sufficient for the considered application, stop, otherwise decrease the value of $\epsilon$. Thereby, the final set of boxes can be reused as a starting decomposition, enlarged by the set of boxes that have a minimal side length in between the current and the former value of $\epsilon$. Even if the value of $\epsilon$ is reduced several times, we expect this variant to %run faster than the one presented in this paper 
be beneficial since the intermediate number of boxes to be stored and handled is smaller. 

Another idea for further speeding up computational time is a parallelization of the algorithm by sending a certain set of boxes to a certain thread. Moreover, one could also use a hybrid sandwich approach which first computes representative points on the convex boundary and then refines in the non-convex parts of the Pareto front if needed. 
Whether these concepts lead to a significant improvement is left for future research. 

\section*{Acknowledgements}
We thank Esther Bonacker for providing us with the radiotherapy example case, which was first discussed in her PhD thesis \cite{Bon19}.

%\section*{References}
%\bibliographystyle{apalike} %abbrvnat
%\begin{verbatim}
\bibliographystyle{tfnlm}
\bibliography{references}
%\end{verbatim}
%\bibliography{interacttfpsample}

\end{document}